\documentclass[a4paper, oneside, 11pt]{article}%
\usepackage{amsmath}
\usepackage{amsfonts}
\usepackage{amssymb}
\usepackage{graphicx}
\usepackage{appendix}
\usepackage[square, numbers, sort&compress]{natbib}
\usepackage{color}%
\providecommand{\U}[1]{\protect\rule{.1in}{.1in}}

\newtheorem{theorem}{Theorem}[section]

\newtheorem{definition}[theorem]{Definition}
\newtheorem{assumption}{Assumption}
\newtheorem{example}[theorem]{Example}

\newtheorem{lemma}[theorem]{Lemma}

\newtheorem{remark}[theorem]{Remark}

\numberwithin{equation}{section}

\newcommand{\E}{{\mathbb E}}
\newcommand{\R}{{\mathbb R}}
\newcommand{\pf}{\noindent\textbf{Proof:} }
\newcommand{\eof}{\hfill{$\Box$}}

\newcommand{\BMO}{L^{2, \;\mathrm{BMO}}_{\mathcal F^{W}}(0, T;\mathbb{R}^{n})}
\usepackage{url}
\usepackage{hyperref}
\hypersetup{
colorlinks=true, 
breaklinks=true, 
urlcolor=green , 
linkcolor=red, 
citecolor=blue, 
pdftitle={}, 
pdfauthor={}, 
pdfsubject={} 
}
\newcommand{\esssup}{\ensuremath{\operatorname*{ess\;sup}}}
\newcommand{\essinf}{\ensuremath{\operatorname*{ess\;inf}}}

\newcommand{\cM}{\ensuremath{\mathcal{M}}}

\newcommand{\argmax}{\ensuremath{\operatorname*{argmax}}}

\newcommand{\red}[1]{{\color{red}#1}}

\usepackage[normalem]{ulem}

\renewcommand{\geq}{\geqslant}
\renewcommand{\leq}{\leqslant}

\begin{document}

\title{Optimal consumption-investment with coupled constraints\\
on consumption and investment strategies\\
in a regime switching market with random coefficients}
\author{Ying Hu \thanks{Univ Rennes, CNRS, IRMAR-UMR 6625, F-35000 Rennes, France. Partially supported by Lebesgue
Center of Mathematics \textquotedblleft Investissements d'avenir\textquotedblright program-ANR-11-LABX-0020-01. Email:
\texttt{ying.hu@univ-rennes1.fr }}
\and Xiaomin Shi\thanks{School of Statistics and Mathematics, Shandong University of Finance and Economics, Jinan
250100, China. Partially supported by NSFC (No.~11801315), NSF of Shandong Province (No.~ZR2018QA001). Email: \texttt{shixm@mail.sdu.edu.cn}}
\and Zuo Quan Xu\thanks{Department of Applied Mathematics, The Hong Kong Polytechnic University, Kowloon, Hong Kong.
Partially supported by NSFC (No.~11971409), Hong Kong
RGC (GRF ~15202421 and 15204622), The PolyU-SDU Joint Research Center on Financial Mathematics, The CAS AMSS-PolyU Joint Laboratory of Applied Mathematics, The Research Centre for Quantitative Finance, The Hong Kong Polytechnic University. Email: \texttt{maxu@polyu.edu.hk}}}
\maketitle

This paper studies finite-time optimal consumption-investment problems with power, logarithmic and exponential utilities, in a regime switching market with random coefficients, subject to coupled constraints on the consumption and investment strategies. We provide explicit optimal consumption-investment strategies and optimal values for the problems in terms of the solutions to some diagonally quadratic backward stochastic differential equation (BSDE) systems and linear BSDE systems with unbound coefficients.
Some of these BSDEs are new in the literature and solving them is one of the main theoretical contributions of this paper. We accomplish the latter by applying the truncation, approximation technique to get some a priori uniformly lower and upper bounds for their solutions.



\addcontentsline{toc}{section}{\hspace*{1.8em}Abstract}

\section{Introduction}
In a consumption-investment model, one aims at maximizing his utility of consumption and wealth by choosing the best consumption and investment strategies in a financial market. Following the pioneering work of Merton \cite{Merton},
a large volume of works has been done on the model. In particular, various constraints such as bankruptcy prohibition, subsistence consumption requirement, wealth-dependent investment and consumption constraints are introduced into the model; see, e.g., Guan, Xu and Yi \cite{GXY}, Sethi \cite{Sethi}, Xu and Yi \cite{XY}, Zariphopoulou \cite{Za}.
All these papers focus on Markovian markets.
The model has also been considered in non-Markovian markets with general random coefficients, by Cox and Huang \cite{CHu}, Karatzas, Lehoczky and Shreve \cite{KLS}, Matoussi and Xing \cite{MX18}, Xing \cite{X17}. The method used in \cite{CHu} and \cite{KLS} is known as martingale method nowadays.
Please refer to the seminar monograph Karatzas and Shreve \cite{KS} for a systematic account on this method for more advanced (incomplete/constrained) markets.

On the other hand, Markov chain is usually adopted to reflect, at the macroeconomic level, the market status (such as bull and bear) in the literature.
This is the famous regime switching model of market; see Hamilton \cite{H89}. Utility maximization and mean-variance models in a regime switching market have been studied by many researchers; see, e.g., \cite{Ba, Liu, SC, Shin, WSZ, Za1, ZY,ZY03,YZ04}. In these works, the market parameters are assumed to be deterministic functions of time for each given regime.
This allows the authors to apply ordinary or partial differential equation (ODE or PDE, for short) method to solve their problems.
In practice, however, even in a bull market, the market parameters, such as the interest rate, stock appreciation rates and volatilities are affected by the uncertainties caused by many factors such as politics, economic, legal, military, corporate governance.
Thus, it is too restrictive to set market parameters as deterministic function of time even if the market status is known. It is necessary to allow the market parameters to depend on not only the Markov chain but also other random resources. With this in mind, the authors consider stochastic linear quadratic control and mean-variance problems in \cite{HSX,HSX1,HSX2}.
Along this framework, this paper generalizes the consumption-investment model of Cheridito and Hu  \cite{CH} to a regime switching market with random coefficients.

This paper studies optimal consumption-investment problems with power, logarithmic and exponential utilities in a regime switching market with parameters depending on both a Markov chain and Brownian motion.
Because of the randomness of the market parameters, the ODE and PDE methods will no longer work for our problem.
To deal with the randomness, it is necessary to apply the backward stochastic differential equation (BSDE) theory of Pardoux and Peng \cite{PP}.
Since its inception in 1990, BSDE theory has become a powerful tool in dealing with modern financial engineering problems such as portfolio selection and asset pricing problems.
For instance, Rouge and El Karoui \cite{RE} characterize the price equation via exponential utility maximization and quadratic BSDEs. By only requiring the strategies take their values in some closed sets, Hu, Imkeller and M$\mathrm{\ddot{u}}$ller \cite{HIM} give solutions to exponential/power/logarithmic utility maximization of terminal wealth and Cheridito and Hu \cite{CH} take into account the intermediate consumption.
Becherer \cite{Be} and Morlais \cite{Mo} study exponential utility optimization and indifference valuation with jumps using quadratic BSDEs driven by random measures. Taking ambiguity and time-consistent ambiguity-averse preferences into account, Laeven and Stadye \cite{LS} investigate the robust portfolio choice and indifference valuation by quadratic BSDEs with infinite activity jumps. Kramkov and Pulido \cite{KP} solve a quadratic BSDE system arising from a price impact model.
Because of regime switching, the BSDE systems in our problems are all multidimensional, making them harder to study compared to one dimensional case.

Meanwhile,
in the study of consumption-investment problems, one often assumes that the consumption and investment strategies are subject to separate constraints; see, e.g., \cite{CH}.
In practice, however, the constraints on them may be coupled together. For instance, when there is no individual constraint on one's consumption and investment strategies but he is not allowed to borrow money, then the constraints are coupled together. This renders the study of the coupled constraints on the consumption and investment strategies, which is, to the best of our knowledge, still rare in the literature, partially due to its complexity. In this paper, different from \cite{CH}, we consider general, not necessarily product, coupled constraints on consumption and investment strategies. This brings new mathematical challenges to our analysis for the related BSDE systems.

We provide explicit optimal consumption and investment strategies in terms of solutions to some BSDE systems.
Because of the emergence of regime switching, the BSDE systems in our model are actually diagonally quadratic coupled through the generator of the Markov chain. Different from \cite{CH}, the quadratic BSDE systems in our paper could not be directly covered by existing literature. Solving these systems is the key theoretical contribution of this paper and we accomplish this by first doing some approximation and truncation of the generators, then finding uniformly lower and upper bounds, and eventually taking limit to obtain the desired solutions. The uniqueness of their solutions are also proved by pure BSDE techniques, so the method may be applied other problems in BSDE theory.

The rest part of this paper is organized as follows.
At the end of this section, we introduce some notations and recall some facts about BMO martingales that will be used frequently in the subsequent analysis.
In Section \ref{fm}, we present the financial market and formulate the consumption-investment problem in a regime switching market with random coefficients. In Sections \ref{powU}, \ref{logU}, \ref{expU}, we solve the problem for the power, logarithmic and exponential utilities, respectively. For each utility,
we first specify the consumption-investment constraint and define the admissible strategies.
We then accomplish the solvability of some related BSDE systems. With the aid of their solutions, we finally provide the optimal consumption and investment strategies and the optimal value of the problem. We conclude the paper in Section \ref{conclude}.


\subsection*{Notation}
Let $(\Omega, \mathcal F, \mathbb{P})$ be a fixed complete probability space on which are defined a standard $n$-dimensional Brownian motion $W_t=(W_{1, t}, \ldots, W_{n, t})'$ and a continuous-time stationary Markov chain $\alpha_t$ valued in a finite state space $\mathcal M=\{1, 2, \ldots, \ell\}$ with $\ell\geq 1$. We assume $\{W_t\}_{t\geq0}$ and $\{\alpha_t\}_{t\geq0}$ are independent processes. The Markov chain has a generator $Q=(q^{ij})_{\ell\times \ell}$ with $q^{ij}\geq 0$ for $i\neq j$ and $\sum_{j=1}^{\ell}q^{ij}=0$ for every $i\in\mathcal{M}$.
Define the filtrations $\mathcal F_t=\sigma\{W_s, \alpha_s: 0\leq s\leq t\}\bigvee\mathcal{N}$ and $\mathcal F^W_t=\sigma\{W_s: 0\leq s\leq t\}\bigvee\mathcal{N}$, where $\mathcal{N}$ is the totality of all the $\mathbb{P}$-null sets of $\mathcal{F}$.

Throughout this paper, we denote by $\R^n$ the set of $n$-dimensional column vectors, by $\R^n_+$ the set of vectors in $\R^n$ whose components are nonnegative, by $\R^{m\times n}$ the set of $m\times n$ real matrices, and by $\mathbb{S}^n$ the set of symmetric $n\times n$ real matrices. For $x\in\R$, we define $x^+:=\max\{x, 0\}$, and $x^-:=\max\{-x, 0\}$. If $M=(m_{ij})\in \R^{m\times n}$, we denote its transpose by $M'$, and its norm by $|M|=\sqrt{\sum_{ij}m_{ij}^2}$. If $M\in\mathbb{S}^n$ is positive definite (positive semidefinite) , we write $M>$ ($\geq$) $0.$ We write $A>$ ($\geq$) $B$ if $A, B\in\mathbb{S}^n$ and $A-B>$ ($\geq$) $0.$
We use the following spaces throughout the paper:
\begin{align*}
L^2_{\mathcal{F}_T}(\Omega;\mathbb{R})&=\Big\{\xi:\Omega\rightarrow
\mathbb{R}\;\Big|\;\xi\mbox { is }\mathcal{F}_{T}\mbox{-measurable, and }\E\big(|\xi|^{2}\big)%
<\infty\Big\}, \\
L^{\infty}_{\mathcal{F}_T}(\Omega;\mathbb{R})&=\Big\{\xi:\Omega\rightarrow
\mathbb{R}\;\Big|\;\xi\mbox { is }\mathcal{F}_{T}\mbox{-measurable, and essentially bounded}\Big\}, \\
L^{2}_{\mathcal F}(0, T;\mathbb{R})&=\Big\{\phi:[0, T]\times\Omega\rightarrow
\mathbb{R}\;\Big|\;\phi\mbox{ is an }\{\mathcal{F}%
_{t}\}_{t\geq0}\mbox{-predictable process with }\\
&\qquad\quad\;\E\int_{0}^{T}|\phi_t|^{2}dt<\infty
\Big\}, \\
L^{1}_{\mathcal F}(0, T;\mathbb{R})&=\Big\{\phi:[0, T]\times\Omega\rightarrow
\mathbb{R}\;\Big|\;\phi\mbox{ is an }\{\mathcal{F}%
_{t}\}_{t\geq0}\mbox{-predictable process with }\\
&\qquad\quad\;\E\int_{0}^{T}|\phi_t|dt<\infty
\Big\}, \\
L^{\infty}_{\mathcal{F}}(0, T;\mathbb{R})&=\Big\{\phi:[0, T]\times\Omega
\rightarrow\mathbb{R}\;\Big|\;\phi\mbox{ is an }\{\mathcal{F}%
_{t}\}_{t\geq0}\mbox{-predictable essentially}\\
&\qquad\mbox{\quad \ bounded process} \Big\}.
\end{align*}
These definitions are generalized in the obvious way to the cases that $\mathcal{F}$ is replaced by $\mathcal{F}^W$ and $\mathbb{R}$ by $\mathbb{R}^n$, $\mathbb{R}^{n\times m}$ or $\mathbb{S}^m$.

In our analysis, some arguments such $s$, $t$, $\omega$, as well as some terms including ``almost surely" (a.s.) and ``almost everywhere" (a.e.) may be suppressed for notation simplicity in some circumstances when no confusion occurs. 

\subsection*{BMO martingales}
We recall some facts about BMO martingales (see Kazamaki \cite{Ka}). They will be used in our subsequent analysis.

For a $\Lambda \in L^{2}_{\mathcal F^{W}}(0, T;\mathbb{R}^{n})$, the process $\int_0^\cdot \Lambda_s'dW_s$ is called a BMO martingale (on $[0,T]$) if there exists a constant $K>0$ such that
\[\mathbb{E}\bigg[\int_\tau^T|\Lambda_s|^2ds\;\bigg|\;\mathcal F_\tau^W\bigg]\leq K\]
for all $\{\mathcal{F}_t^W\}_{t\geq 0}$-stopping times $\tau\leq T$. The set of BMO martingales is defined as
\begin{align*}
\BMO&=\bigg\{\Lambda \in L^{2}_{\mathcal F^{W}}(0, T;\mathbb{R}^{n}) \;\bigg|\; \int_0^\cdot\Lambda_s'dW_s \mbox{ is a BMO martingale on $[0, T]$}\bigg\}.
\end{align*}
The following two important properties of BMO martingales will be used in our below arguments without claim.
\begin{itemize}
\item If $\int_0^\cdot \Lambda_s'dW_s$ is a BMO martingale, then Dol$\acute{\mathrm{e}}$ans-Dade stochastic exponential
$$\mathcal E\bigg(\int_0^\cdot \Lambda_s'dW_s\bigg)$$
 is a uniformly integrable martingale.
 \item
If $\int_0^\cdot \Lambda_s'dW_s$ and $\int_0^\cdot Z_s'dW_s$ are both BMO martingales, then $\widetilde W:=W-\int_0^\cdot Z_sds$ is a standard Brownian motion, and $\int_0^\cdot \Lambda_s'd\widetilde W_s$ is a BMO martingale under the probability measure $\widetilde{\mathbb{P}}$ defined by $$\frac{d\widetilde{\mathbb{P}}}{d\mathbb{P}}\bigg|_{\mathcal{F}^{W}_T}=\mathcal E \bigg(\int_0^T Z_s'dW_s\bigg).$$
\end{itemize}

\section{Problem formulation}\label{fm}

\subsection{Financial market}
Consider a financial market consisting of a risk-free asset (the money market
instrument or bond) whose price is $S_{0}$ and $m$ risky securities (the
stocks) whose prices are $S_{1}, \ldots, S_{m}$. Assume $m\leq n$, i.e., the number of risky securities is no more than
that of the market risk resources (namely the Brownian motion). When $m<n$ the market is incomplete.
These asset prices are driven by SDEs:
\begin{align*}
\begin{cases}
dS_{0, t}=r_{t}^{\alpha_t}S_{0, t}dt, \\
S_{0, 0}=s_0,
\end{cases}
\end{align*}
and
\begin{align*}
\begin{cases}
dS_{k, t}=S_{k, t}\Big(\mu_{k, t}^{ \alpha_t}dt+\sum\limits_{j=1}^n\sigma_{kj, t}^{ \alpha_t}dW_{j, t}\Big), \\
S_{k, 0}=s_k,
\end{cases}
\end{align*}
where $r^i_t$ is the interest rate process and $\mu_{k, t}^i$ and $\sigma_{k, t}^i:=(\sigma_{k1, t}^i, \ldots, \sigma_{kn, t}^i)$ are the appreciation rate process and volatility rate process of the $k$th risky security corresponding to a market regime $\alpha_t=i$, for every $k=1, \ldots, m$ and $i\in\cM$. Recall that $\alpha_t$ follows a continuous-time stationary Markov chain valued in a finite state space $\mathcal M=\{1, 2, \ldots, \ell\}$ with $\ell\geq 1$. When $\ell=1$, there is no regime switching and the market becomes the classical Black-Scholes market. In our below argument we assume $\ell>1$, although all the results remain true if $\ell=1$.

Define the appreciate vector
\begin{align*}
\mu_{t}^i=(\mu_{1, t}^i, \ldots, \mu_{m, t}^i)',
\end{align*}
and volatility matrix
\begin{align*}
\sigma_{t}^i=
\left(
\begin{array}{c}
\sigma_{1, t}^i\\
\vdots\\
\sigma_{m, t}^i\\
\end{array}
\right)
\equiv (\sigma_{kj, t}^i)_{m\times n}, \ \text{for}\ \text{each} \quad i\in\cM.
\end{align*}

\subsection{Optimal investment-consumption problem}
Consider a small investor, whose actions cannot affect the asset prices. He will decide at every time
$t\in[0, T]$ the \emph{proportion} $\pi_{j, t}$ of his wealth to invest in the $j$th risky asset, $j=1, \ldots, m$, as well as the \emph{proportion} $ c_t$ of his wealth to consume.
The vector process $\pi_{t}:=(\pi_{1, t}, \ldots, \pi_{m, t})'$ is called a portfolio of the investor.
The pair $(\pi,c)$ is called a consumption-investment strategy.
The investor's self-financing wealth process $X$ corresponding to a consumption-investment strategy $(\pi,c)$ is a strong solution of the SDE (see, e.g., \cite{KS}):
\begin{align}
\label{wealth}
\begin{cases}
dX_t=X_t[r_{t}^{\alpha_t}+\pi_t'b_{t}^{\alpha_t}- c_t ]dt+X_t\pi_t'\sigma_{t}^{\alpha_t}dW_t, \\
X_0=x>0, \ \alpha_0=i_0\in\cM,
\end{cases}
\end{align}
where $b_{t}^{\alpha_t}:=\mu_{t}^{\alpha_t}-r_{t}^{\alpha_t}\mathbf{1}_{m}$ and $\mathbf{1}_{m}$ is the $m$-dimensional vector with all entries being one. By It\^{o}'s lemma, one can easily show that the process $X$ is always positive, hence the
investor would never be bankrupt.

The investor's problem is to maximize
\begin{align}
J(x, i_0;\pi, c)&:=\E\bigg[\int_0^T e^{-\int_0^t\rho_{s}^{\alpha_s}ds}U(c_tX_t)dt+e^{-\int_0^T\rho_{s}^{\alpha_s}ds}U(X_T)\bigg], \quad \mathrm{ s.t.} \;\; (\pi, c)\in \mathcal{U},
\label{optm}%
\end{align}
and determine the value function $$V(x, i_0):=\sup_{(\pi, c)\in\mathcal{U}}J(x, i_0;\pi, c), $$ where $\rho^i\in L^\infty_{\mathcal F^W}(0, T;\mathbb R)$, $i\in\cM$, are the discount factor processes, $U$ is the utility function of the investor, and $\mathcal{U}$ is the admissible set of consumption-investment strategies. The utility $U$ and admissible set $\mathcal{U}$ will be defined in the sequel case by case. In particular, the consumption-investment strategies can be subject to fairly general constraints.
Let $\Theta$ denote the constraint set for them, which is assumed to be a given closed nonempty set in $\mathbb{R}^m\times\mathbb{R}_{+}$.

\begin{example}
Here are some important and interesting examples for the constraint set $\Theta$ in financial practice.
\begin{itemize}
\item
If the investor are not allowed to borrow money, then the total consumption and investment cannot beyond 100\%, that is,
\begin{align*}
\Theta&=\Big\{(\pi, c)\in\mathbb{R}^m\times\mathbb{R}_{+}\;\Big|\; \sum_{j=1}^{m}\pi_j+c\leq 1\Big\}.
\end{align*}
Because the constraints on the consumption and investment strategies are coupled together, this will bring some mathematical challenges to our analysis for the related BSDE systems below.

\item If no shorting is allowed in the stocks, then
\begin{align*}
\Theta &=\Big\{(\pi, c)\in\mathbb{R}^m\times\mathbb{R}_{+}\;\Big|\; \pi_{j}\geq 0, \;j=1, \ldots, m\Big\}.
\end{align*}
In this case $\Theta$ is cone.
\item There may be restrictions on the investment strategies on some stocks; for instance
\begin{align*}
\Theta &=\Big\{(\pi, c)\in\mathbb{R}^m\times\mathbb{R}_{+}\;\Big|\; \pi_{j}\in[d_{j}, e_{j}], \;j=1, \ldots, m\Big\};
\end{align*}
or constraint on the consumption strategies such as
\begin{align*}
\Theta &=\Big\{(\pi, c)\in\mathbb{R}^m\times\mathbb{R}_{+}\;\Big|\; c\in[0, \tfrac{1}{5}]\Big\};
\end{align*}
or coupled constraints on both the investment and consumption strategies
\begin{align*}
\Theta &=\Big\{(\pi, c)\in\mathbb{R}^m\times\mathbb{R}_{+}\;\Big|\; \; c\in[0, \tfrac{1}{5}],\;\pi_{j}\in[d_{j}, e_{j}], \;j=1, \ldots, m\Big\}.
\end{align*}
\end{itemize}
\end{example}
We put the following standard assumption for the market parameters.
\begin{assumption} \label{assu1}
For all $i\in\cM$,
\begin{align*}
\begin{cases}
r^{i}\in L_{\mathbb{F}^W}^\infty(0, T;\mathbb{R}), \
\mu^{i}\in L_{\mathbb{F}^W}^\infty(0, T;\mathbb{R}^m), \\
\sigma^{i}\in L_{\mathbb{F}^W}^\infty(0, T;\mathbb{R}^{m\times n}), \
\rho^i\in L_{\mathbb{F}^W}^\infty(0, T;\mathbb{R}).
\end{cases}
\end{align*}
Also, there exists a constant $\delta>0$ such that $\sigma^{i}(\sigma^{i})'\geq\delta I_{m}$ for all $i\in\cM$, where $I_m$ denotes the $m$-dimensional identity matrix.
\end{assumption}

Due to different features of the power, logarithmic and exponential utilities, they must be dealt with different methods.
We study them in the subsequent Sections \ref{powU}, \ref{logU} and \ref{expU}, respectively.

\section{Power utility}\label{powU}
In this section, we assume that the investor's utility function is
\[U(x)=\frac{1}{\gamma}x^\gamma,\quad x>0,\quad \gamma\in(-\infty, 0)\cup(0, 1)\]
and
the admissible consumption-investment set is defined as
\begin{align*}
\mathcal U=\bigg\{(\pi, c)\;\bigg|\; \int_0^T (|\pi_{t}|^2+c_{t})dt<\infty, a.s. , \ (\pi_t(\omega), c_t(\omega))\in\Theta, \ a.e. , \ a.s, \bigg\},
\end{align*}
if $\gamma\in(0, 1)$; and
\begin{align*}
\mathcal U&=\bigg\{(\pi, c)\;\bigg|\; \int_0^T (|\pi_{t}|^2+c_{t})dt<\infty, a.s., \ (\pi_t(\omega), c_t(\omega))\in\Theta, \ a.e. , \ a.s, \\
&\qquad\qquad \qquad (X_{t}^{\gamma})_{0\leq t\leq T} \ \mbox{belongs to class (D) and } \ \E\bigg[\int_0^T (c_{t}X_{t})^{\gamma}dt\bigg]<\infty\bigg\},
\end{align*}
if $\gamma<0$.\footnote{We say $(X_{t}^{\gamma})_{0\leq t\leq T}$ belongs to class (D) if $\{X_{\tau}^{\gamma}:\tau \ \mbox{stopping time valued in } \ [0,T]\}$ is a family of uniformly integrable random variables.}

\begin{assumption} \label{assu0}
If $\gamma\in(0,1)$, then no consumption or investment is always permitted, namely $(0\mathbf{1}_{m}, 0)\in\Theta$.
If $\gamma<0$, then positive consumption is always needed (for otherwise the admissible set $\mathcal U$ ie empty),
there exists $\varepsilon>0$ such that $(0\mathbf{1}_{m}, \varepsilon)\in\Theta$.
\end{assumption}

\bigskip

To tackle problem \eqref{optm}, we introduce the following $\ell$-dimensional BSDE system:
\begin{align}
\label{P1}
\begin{cases}
dP_t^i=-[ f^i(P_t^i, \Lambda_t^i)-(\rho^i-\gamma r^i)P^i+\sum_{j=1}^\ell q^{ij}P^j]dt+(\Lambda_t^i)'dW, \\
 P_T^i=1, \\
 P^i>0, \quad i\in\cM,
\end{cases}
\end{align}
where, for any $(P, \Lambda)\in\mathbb{R}_+\times\mathbb{R}^n$,
\begin{align*}
f^i(P, \Lambda)&=\gamma\esssup_{(\pi, c)\in\Theta}\bigg[-\frac{1-\gamma}{2}|\pi'\sigma^i|^2 P+\pi'(Pb^i+\sigma^i\Lambda)+\frac{c^\gamma}{\gamma}-Pc\bigg].
\end{align*}

\begin{remark}
If there is no consumption-investment constraint, namely $\Theta=\mathbb{R}^n\times\mathbb{R}_+$, then
\begin{align*}
 f^i(P, \Lambda)
&=\frac{1}{2}\frac{\gamma}{1-\gamma}\frac{1}{P}(Pb+\sigma\Lambda)'(\sigma\sigma')^{-1}(Pb+\sigma\Lambda)+(1-\gamma)P^{-\frac{\gamma}{1-\gamma}}.
\end{align*}
\end{remark}

\begin{definition}\label{def}
A vector process $(P^i, \Lambda^i)_{i\in\cM}$ is called a solution to the $\ell$-dimensional BSDE \eqref{P1}, if it satisfies \eqref{P1}, and $(P^i, \Lambda^i)\in L^\infty_{\mathcal{F}^W}(0, T; \mathbb{R})\times \BMO$ for all $i\in\cM$. A solution $(P^i, \Lambda^i)_{i\in\cM}$ to \eqref{P1} is called uniformly positive if $P^i_t\geq \delta$, for all $t\in[0, T]$, a.s. with some deterministic constant $\delta>0$.
\end{definition}


The following comparison theorem for multidimensional BSDE systems firstly appeared in \cite{HP} (one can find a concise version in \cite[Lemma 2.2]{HLT} or \cite[Lemma 3.4]{HSX}). We shall use it frequently in the study of BSDEs.
\begin{lemma}
\label{comparison}
Suppose $(Y^i, Z^i)_{i\in\cM}$ and $(\overline Y^i, \overline Z^i)_{i\in\cM}$ satisfy the following two $\ell$-dimensional BSDE systems, respectively:
\begin{align*}
Y^i_t=\xi^i+\int_t^T g^i(s, Y^i_s, Y^{-i}_s, Z^i_s)ds-\int_t^T (Z^i_s)'dW_s, \ \mbox{ for all $i\in\cM$;}
\end{align*}
and
\begin{align*}
\overline Y^i_t=\overline\xi^i+\int_t^T \overline g^i(s, \overline Y^{i}_s, \overline Y^{-i}_s, \overline Z^{i}_s)ds-\int_t^T (\overline Z^{i}_s)'dW_s, \ \mbox{ for all $i\in\cM$, }
\end{align*}
where $Y^{-i}_s=(Y^1_s, \ldots, Y^{i-1}_s, Y^{i+1}_s, \ldots, Y^{\ell}_s)$.
Also suppose that, for all $i\in\cM$,
\begin{enumerate}
\item $\xi^i, \ \overline\xi^i\in L^2_{\mathcal{F}^W}(\Omega;\mathbb{R})$, and $\xi^i\leq\overline\xi^i$;
\item there exists a constant $K>0$ such that
\[|g^i(s, y, z)-g^i(s, \overline y, \overline z)|\leq K(|y-\overline y|+|z-\overline z|), \]
for any $z, \overline z\in\mathbb{R}^n$, $y=(y^i, y^{-i})$, $\overline y=(\overline y^i, \overline y^{-i})\in\mathbb{R}^{\ell}$;

\item $g^i(s, y, z)$ is nondecreasing in every $y^j$, $j\neq i \in\cM$; and
\item $ g^i(s, \overline Y^{i}_s, \overline Y^{-i}_s, \overline Z^{i}_s)\leq \overline g^i(s, \overline Y^{i}_s, \overline Y^{-i}_s, \overline Z^{i}_s)$.
\end{enumerate}
Then $Y^{i}_t\leq \overline Y^{i}_t$ a.s. for all $t\in[0, T]$ and $i\in\cM$.
\end{lemma}
\begin{remark}
If conditions 2 and 3 hold for $\overline g^i$ instead of $g^i$, and condition 4 is replaced by $$g^i(s, Y^{i}_s, Y^{-i}_s, Z^{i}_s)\leq \overline g^i(s, Y^{i}_s, Y^{-i}_s, Z^{i}_s)$$ in the above lemma, then the conclusion still holds.
\end{remark}

\begin{theorem}\label{Th:P}
Suppose $\gamma\in(0, 1)$ and Assumptions \ref{assu1} and \ref{assu0} hold.
Then there is a unique uniformly positive solution $(P^i, \Lambda^i)_{i\in\cM}$ to the BSDE \eqref{P1}.
\end{theorem}
\pf
Let $a>1$ be a large constant such that
\begin{align}\label{a1}
-\rho^i+\gamma r^i\geq -a,
\end{align}
and
\begin{align}\label{a2}
\frac{\gamma}{2(1-\gamma)}(b^i)'(\sigma^i(\sigma^i)')^{-1}b^i-\rho^i+\gamma r^i\leq a(1-\gamma),
\end{align}
hold simultaneously for all $i\in\cM$.
Denote $a_1:=e^{-aT}$ and $a_2:=2e^{aT}$.

Let $g:\R \rightarrow[0, 1]$ be a smooth truncation function satisfying $g(x)=0$ for $x\in(-\infty, \frac{1}{2}a_1]$, and $g(x)=1$ for $x\in[a_1, +\infty)$.
For $k\geq1$, $(t, P, \Lambda)\in[0, T]\times\mathbb R\times\mathbb{R}^n$, $ i\in\cM$, define
\[
f^{k, i}(t, P, \Lambda)=\essinf_{\tilde P\in\mathbb R, \tilde\Lambda\in\mathbb R^n}\Big[f^i (t, \tilde P, \tilde\Lambda)g(\tilde P)+k|P-\tilde P|+k|\Lambda-\tilde\Lambda|\Big].
\]
Then it is uniformly Lipschitz continuous in $(P, \Lambda)$, and increasingly approaches to $f^i(t, P, \Lambda)g(P)$ as $k$ goes to infinity. According to Assumption \ref{assu0}, we have $f^i\geq 0$, hence $f^{k, i}\geq 0$.

The following BSDE system
\begin{align*}
\begin{cases}
dP^{k, i}=-\Big[f^{k, i}(P^{k, i}, \Lambda^{k, i})-(\rho^i-\gamma r^i)P^{k, i}+\sum_{j=1}^\ell q^{ij}P^{k, j}\Big]dt+(\Lambda^{k, i})^\top dW, \\
P^{k, i}_T=1, \ \mbox{ for all $i\in\cM$, }
\end{cases}
\end{align*}
has a Lipschitz generator\footnote{As for BSDE $Y_t=\xi+\int_t^T f(s, Y_s, Z_s)ds-\int_t^T Z_s'dW_s$, $f$ is called the generator and $\xi$ is called the terminal value.}, so it admits a unique solution, denoted by $\big(P^{k, i}, \Lambda^{k, i}\big)_{i\in\cM}$.

It is direct to verify that
\begin{align*}
(\underline P^i_t, \underline \Lambda^i_t)=\Big(e^{-a(T-t)}, 0\Big), \quad i\in\cM
\end{align*}
is the unique solution to the following linear BSDE system:
\begin{align*}
\begin{cases}
d\underline P^i=-\Big[-a\underline P^i+\sum_{j=1}^\ell q^{ij}\underline P^j\Big]dt+(\underline\Lambda^i)'dW, \\
\underline P^i_T=1, \quad i\in\cM.
\end{cases}
\end{align*}
Notice that $f^{k, i}(P, \Lambda)\geq0$ and recall \eqref{a1}, so
\begin{align*}
f^{k, i}(\underline P^{i}, \underline \Lambda^{i})-(\rho^i-\gamma r^i)\underline P^{i}+\sum_{j=1}^\ell q^{ij}\underline P^{j}
&\geq -(\rho^i-\gamma r^i)\underline P^{i}+\sum_{j=1}^\ell q^{ij}\underline P^{j} \geq -a\underline P^{i}+\sum_{j=1}^\ell q^{ij}\underline P^{j}.
\end{align*}
By Lemma \ref{comparison}, we have
\[
a_1\leq e^{-a(T-t)}=\underline P^i_t\leq P^{k, i}_t, \quad i\in\cM.
\]

On the other hand,
\begin{align*}
(\bar P^i, \bar \Lambda^i)=\bigg(\Big(e^{a(T-t)}+\frac{1}{a}(e^{a(T-t)}-1)\Big)^{1-\gamma}, \; 0\bigg), \quad i\in\cM,
\end{align*}
is a solution to the following BSDE system:
\begin{align*}
\begin{cases}
d\bar P^i=-\Big[ a(1-\gamma) \bar P^i +(1-\gamma)(\bar P^i)^{-\frac{\gamma}{1-\gamma}}+\sum_{j=1}^\ell q^{ij}\bar P^j\Big]dt+(\bar\Lambda^i)'dW, \\
\bar P^i_T=1, \quad i\in\cM.
\end{cases}
\end{align*}
Notice, for $P> 0$,
\begin{align}\label{quadratic}
f^{k, i}(t, P, \Lambda)&\leq f^{i}(t, P, \Lambda)g(P)\nonumber\\
&\leq \gamma \esssup_{(\pi, c)\in\R^n\times\R_+}\Big[-\frac{1-\gamma}{2}|\pi'\sigma^i|^2 P+\pi'(Pb^i+\sigma^i\Lambda)+\frac{c^\gamma}{\gamma}-Pc\Big]g(P)\nonumber\\
&=\frac{1}{2}\frac{\gamma}{1-\gamma}\frac{1}{P}(Pb+\sigma\Lambda)'(\sigma\sigma')^{-1}(Pb+\sigma\Lambda)g(P)
+(1-\gamma)P^{-\frac{\gamma}{1-\gamma}}g(P)\nonumber\\
&\leq \frac{\gamma}{2(1-\gamma)}(b^i)'(\sigma^i(\sigma^i)')^{-1}b^i P +\frac{\gamma}{1-\gamma}(b^i)'(\sigma^i(\sigma^i)')^{-1}\sigma^i\Lambda\nonumber\\
&\qquad+\frac{\gamma}{2(1-\gamma)}\frac{1}{P}(\sigma^i\Lambda)'(\sigma^i(\sigma^i)')^{-1}\sigma^i\Lambda
+(1-\gamma)P^{-\frac{\gamma}{1-\gamma}}.
\end{align}
Hence, by \eqref{a2},
\begin{align*}
&\quad f^{k, i}(\bar P^{i}, \bar \Lambda^{i})-(\rho^i-\gamma r^i)\bar P^{i}+\sum_{j=1}^\ell q^{ij}\bar P^{j}\\
&\leq \frac{\gamma}{2(1-\gamma)}(b^i)'(\sigma^i(\sigma^i)')^{-1}b^i \bar P^i
+(1-\gamma)(\bar P^i)^{-\frac{\gamma}{1-\gamma}}-(\rho^i-\gamma r^i)\bar P^{i}+\sum_{j=1}^\ell q^{ij}\bar P^{j}\\
&\leq a(1-\gamma)\bar P^i+(1-\gamma)(\bar P^i)^{-\frac{\gamma}{1-\gamma}}+\sum_{j=1}^\ell q^{ij}\bar P^{j}.
\end{align*}
Then by Lemma \ref{comparison} again, we have
\[
 P^{k, i}_t\leq \bar P^i_t\leq 2e^{a T}=a_2,
\]
and $P^{k, i}$ is increasing in $k$, for each $i\in\cM$.

Let $P^i_t=\lim\limits_{k\rightarrow\infty}P^{k, i}_t$, $i\in\cM$. Since $a_1$ and $a_2$ are independent of $k$, $a_1\leq P^i_t\leq a_2$. Recalling $f^{k, i}\geq0$, \eqref{quadratic} and the role of the truncation function $g$, we can regard $\big(P^{k, i}, \Lambda^{k, i}\big)$ as the solution of a scalar-valued quadratic BSDE for each $i\in\cM$. Thus by \cite[Proposition 2.4]{Ko}, there exists a process $\Lambda\in L^{2}_{\mathcal F^{W}}(0, T;\mathbb{R}^{n\times \ell})$ such that $(P, \Lambda)$ satisfies the BSDE \eqref{P1}.

Applying It\^{o}'s formula to $(P^i-a_2)^2$, we obtain, for any stopping times $\tau\leq T$,
\begin{align*}
\E\bigg[\int_\tau^T|\Lambda^i|^2\;\bigg|\;\mathcal{F}_\tau^W\bigg]&=(1-a_2)^2-(P_\tau^i-a_2)^2\\
&\qquad+2\E \bigg[\int_{\tau}^T(P^i-a_2)\Big[f^i(P^i, \Lambda^i)-(\rho^i-\gamma r^i)P^i+\sum_{j=1}^\ell q^{ij}P^j\Big]ds \;\bigg|\;\mathcal{F}_\tau^W\bigg]\\
&\leq (1-a_2)^2-(P_\tau^i-a_2)^2\\
&\qquad+2\E\bigg[ \int_{\tau}^T(P^i-a_2)\Big[-(\rho^i-\gamma r^i)P^i+\sum_{j=1}^\ell q^{ij}P^j\Big]ds \;\bigg|\;\mathcal{F}_\tau^W\bigg]\\
&\leq (1-a_2)^2+2\E\bigg[\int_{0}^T|P^i-a_2|\Big|-(\rho^i-\gamma r^i)P^i+\sum_{j=1}^\ell q^{ij}P^j\Big|ds \;\bigg|\;\mathcal{F}_\tau^W\bigg],
\end{align*}
where we used the fact that $P^i\leq a_2$ and $f^i\geq 0$ to get the first inequity.
Because $a_1\leq P^i\leq a_2$ and Assumption \ref{assu1}, we see the right hand side is upper bounded by a constant.
Hence $\Lambda^i\in\BMO$, for all $i\in\cM$.
We have now established the existence of the solution.

Next, let us prove the uniqueness.
Suppose $(P^i, \ \Lambda^i)_{i\in\cM}$, $(\tilde P^i, \ \tilde\Lambda^i)_{i\in\cM}$ are two uniformly positive solutions of $\eqref{P1}$.
For every $i\in\cM$, define processes
\begin{align}\label{logtrans}
(Y_t^i, Z_t^i)=\left(\ln P_t^i, \;\frac{\Lambda_t^i}{P_t^i}\right), \
(\tilde Y_t^i, \tilde Z_t^i)=\left(\ln \tilde P_t^i, \; \frac{\tilde \Lambda_t^i}{\tilde P_t^i}\right), \ \mbox{ for } t\in[0, T].
\end{align}
Then $(Y^i, Z^i)$, $(\tilde Y^i, \tilde Z^i)\in L^\infty_{\mathcal{F}^W}(0, T; \mathbb {R})\times \BMO$, for all $ i\in\cM$. Furthermore, by It\^{o}'s formula, $(Y^i, Z^i)_{i\in\cM}$ and $(\tilde Y^i, \tilde Z^i)_{i\in\cM}$ satisfy the following BSDE system:
\begin{align}\label{Y}
\begin{cases}
dY^i=-\Big[F^i(Y^i, Z^i)+\frac{1}{2}|Z^i|^2-\rho^i+\gamma r^i+\sum_{j=1}^\ell q^{ij}e^{Y^j-Y^i}\Big]dt+(Z^i)'dW, \\
Y^i_T=0, \quad i\in\cM.
\end{cases}
\end{align}
where for any $(Y, Z)\in\mathbb{R}\times\mathbb{R}^n$,
\begin{align*}
F^i(Y, Z)&=\gamma\esssup\limits_{(\pi, c)\in\Theta}\Big[-\frac{1-\gamma}{2}|\pi'\sigma^i|^2 +\pi'(b^i+\sigma^i Z)+\frac{c^\gamma}{\gamma}e^{-Y}-c\Big].
\end{align*}

Because $\ln a_1\leq Y^i\leq \ln a_2$ and use Assumption \ref{assu1}, for any $(\pi, c)\in\Theta$, 
\begin{align*}
&\quad -\frac{1-\gamma }{2}|\pi'\sigma^i|^2 +\pi'(b^i+\sigma^i Z^i)+\frac{c^\gamma}{\gamma}e^{-Y^i}-c\\
&\leq -\delta \frac{1-\gamma }{2}|\pi|^2+K_{1}|\pi|(1+|Z|) +\sup_{c\in\R_+}\Big(\frac{c^\gamma}{\gamma}\frac{1}{a_{1}}-c\Big)\\
&= -\delta \frac{1-\gamma }{2}|\pi|^2+K_{1}|\pi|(1+|Z|) +\frac{1-\gamma}{\gamma}a_{1}^{\frac{1}{1-\gamma}}\\
&<0,
\end{align*}
if $|\pi|>K(1+|Z^i|)$ for sufficient large $K$.\footnote{ Hereafter, we shall use $K$ to represent a generic positive constant independent of $i, \ m, \ n$ and $t$, which can be different from line to line.}
Since $0\mathbf{1}_{m+1}\in\Theta$, $F^i(Y, Z)\geq 0$. Therefore
\begin{align}\label{Fboundeddomain1}
F^i(Y^i, Z^i)&=\gamma\esssup_{\substack{(\pi, c)\in\Theta\\|\pi|\leq K(1+|Z^i|)}}\Big[-\frac{1-\gamma}{2}|\pi'\sigma^i|^2 +\pi'(b^i+\sigma^i Z^i)+\frac{c^\gamma}{\gamma}e^{-Y^i}-c\Big].
\end{align}
This still holds if we replace $K(1+|Z^i|)$ by any bigger value.
A similar expression holds for $F^i(\tilde Y^i, \tilde Z^i)$.

Set $\overline Y^i=Y^i-\tilde Y^i, \ \overline Z^i=Z^i-\tilde Z^i$, for $i\in\cM$.
From \eqref{Fboundeddomain1}, we have
\begin{align*}
\Big|F^i(Y^i, Z^i)-F^i( Y^i, \tilde Z^i)\Big|
&\leq \gamma \esssup_{\substack{(\pi, c)\in\Theta\\|\pi|\leq K(1+|Z^i|+|\tilde Z^i|)}}|\pi'\sigma^i \overline Z^i| \leq K(1+|Z^i|+|\tilde Z^i|) |\overline Z^i|.
\end{align*}
Since $Z^i, \tilde Z^i\in\BMO$, we can define $\beta^i\in\BMO$ in an obvious way such that
$$F^i(Y^i, Z^i)-F^i( Y^i, \tilde Z^i)=(\beta^i)'\overline Z^i,$$
and $$|\beta^i| \leq K(1+|Z^i|+|\tilde Z^i|).$$

Now applying It\^{o}'s formula to $(\overline Y^i)^2$, we deduce that
\begin{align}\label{overlineY}
(\overline Y^i_t)^2&=\int_t^T \Big\{2\overline Y^i\Big[F^i(Y^i, Z^i)-F^i(\tilde Y^i, \tilde Z^i)+\frac{1}{2}(|Z^i|^2-|\tilde Z^i|^2)\nonumber\\
&\qquad+\sum_{j=1}^\ell q^{ij}\Big(e^{Y^j-Y^i}-e^{\tilde Y^j-\tilde Y^i}\Big)\Big]-|\overline Z^i|^2\Big\}ds-\int_t^T2\overline Y^i(\overline Z^i)'dW.
\end{align}
Notice that the map
\[
Y\mapsto F^i(Y, Z), \ Y\in\R_+,
\]
is non-increasing for every $Z\in\R^n$ and $i\in\cM$. Therefore,
\begin{align*} \overline Y^i\big[F^i(Y^i, Z^i)-F^i(\tilde Y^i, \tilde Z^i)\big]
&=\overline Y^i\big[F^i(Y^i, Z^i)-F^i( Y^i, \tilde Z^i)+F^i(Y^i, \tilde Z^i)-F^i(\tilde Y^i, \tilde Z^i)\big]\\
&\leq \overline Y^i\big[F^i(Y^i, Z^i)-F^i( Y^i, \tilde Z^i)\big]\\
&=\overline Y^i(\beta^i)'\overline Z^i.
\end{align*}

For each fixed $i\in\cM$, let us introduce the process
\[
N^i_t=\mathcal E\Big(\int_0^t(\beta^i_s+\frac{1}{2}(Z^i+\tilde Z^i))'dW_s\Big).
\]
Then $N_t^i$ is a uniformly integrable martingale. Notice $Y^i$, $i\in\cM$ are bounded, from \eqref{overlineY},
\begin{align*}
(\overline Y^i_t)^2&\leq \int_t^T \Big\{2\overline Y^i\Big[(\beta^i)'\overline Z^i+\frac{1}{2}(Z^i+\tilde Z^i)'\overline Z^i+K\sum_{j\neq i}^\ell \overline Y^j\Big]\Big\}ds-\int_t^T2\overline Y^i(\overline Z^i)'dW\\
&= \int_t^T 2K\overline Y^i\sum_{j\neq i}^\ell \overline Y^jds-\int_t^T2\overline Y^i(\overline Z^i)'d\widetilde W^i
\end{align*}
where
\[\widetilde W^i_t:=W_t-\int_0^t\Big(\beta^i_s+\frac{1}{2}(Z^i+\tilde Z^i)\Big)ds, \]
is a Brownian motion under the probability $\widetilde{\mathbb P}^{i}$ defined by
\begin{align*}
\frac{d\widetilde{\mathbb P^{i}}}{d\mathbb P}\bigg|_{\mathcal{F}_T^W}=N^i_T.
\end{align*}
Taking expectation $\widetilde\E^i$ w.r.t. the probability measure $\widetilde{\mathbb P}^i$,
\begin{align*}
(\overline Y^i_t)^2&\leq \widetilde\E^i\bigg[\int_t^T2K\overline Y^i\sum_{j\neq i}^\ell \overline Y^jds\;\bigg|\;\mathcal{F}^W_t\bigg] \leq K \widetilde\E^i\bigg[\int_t^T\sum_{j=1}^\ell (\overline Y^j)^2ds\;\bigg|\;\mathcal{F}^W_t\bigg] \leq K\int_t^T \sum_{j=1}^\ell E_s^j ds,
\end{align*}
by the arithmetic-mean and geometric-mean inequality (AM-GM inequality),
where
\[
E^i_t=\underset{\omega\in\Omega}{\esssup} \ { (\overline Y_t^i)^2}.
\]
Taking essential supreme on both sides, we deduce
\[
E_t^i\leq K \int_t^T \sum_{j=1}^{\ell}E_s^jds.
\]
Thus
\[
0\leq \sum_{j=1}^{\ell}E_t^j\leq K\ell\int_t^T \sum_{j=1}^{\ell}E_s^j ds.
\]
We infer from Gronwall's inequality that $\sum_{j=1}^{\ell}E_t^j=0$, so
$\overline Y_t^i=0$ a.s. for $t\in[0, T]$ and all $i\in\cM$. This completes the proof of the uniqueness.
\eof

\begin{remark}
Please note that the solvability of \eqref{P1} and \eqref{Y} {cannot} be directly covered by Fan, Hu and Tang \cite{FHT} since the generators violate the locally Lipschitz condition required in \cite{FHT}.
\end{remark}

\begin{theorem}\label{Th:Pgamma}
Suppose $\gamma<0$ and Assumption \ref{assu1} and \ref{assu0} hold. Then there is a unique uniformly positive solution $(P^i, \Lambda^i)_{i\in\cM}$ to the BSDE \eqref{P1}.
\end{theorem}
\pf
The proof is similar to the procedure of Theorem \ref{Th:P} with different uniformly lower bounds and upper bounds. Hence we will only present how to find these bounds. Other details are left to the interested readers.

Let $a'>0$ be a large constant such that
\begin{align}\label{aprime1}
\frac{\gamma}{2(1-\gamma)}(b^i)'(\sigma^i(\sigma^i)')^{-1}b^i-\rho^i+\gamma r^i\geq -a'(1-\gamma),
\end{align}
and
\begin{align}\label{aprime2}
-\rho^i+\gamma r^i-\gamma \varepsilon\leq a',
\end{align}
hold simultaneously for all $i\in\cM$, where $\varepsilon>0$ is given in Assumption \ref{assu0}.
Denote $a_1':=e^{-a'(1-\gamma)T}$ and $a_2':=e^{a'T}+\frac{\varepsilon^\gamma}{a'}(e^{a'T}-1)$.

Let $g:\R \rightarrow[0, 1]$ be a smooth truncation function satisfying $g(x)=0$ for $x\in (-\infty, \frac{1}{2}a'_1]$, and $g(x)=1$ for $x\in[a'_1, +\infty)$. Notice we have $|xg(x)|\leq |x|$ for all $x\in\R$.

For integer $k\geq -\gamma\varepsilon$, $(t, P, \Lambda)\in[0, T]\times\mathbb R\times\mathbb{R}^n$, $ i\in\cM$, define
\[
f^{k, i}(t, P, \Lambda)=\esssup_{\tilde P\in\mathbb R, \tilde\Lambda\in\mathbb R^n}\Big[f^i (t, \tilde P, \tilde\Lambda)g(\tilde P)-k|P-\tilde P|-k|\Lambda-\tilde\Lambda|\Big].
\]
And consider the following $\ell$-dimensional BSDE:
\begin{align*}
\begin{cases}
dP^{k, i}=-\Big[f^{k, i}(P^{k, i}, \Lambda^{k, i})-(\rho^i-\gamma r^i)P^{k, i}g(P^{k, i})+\sum_{j=1}^\ell q^{ij}P^{k, j}\Big]dt+(\Lambda^{k, i})^\top dW, \\
P^{k, i}_T=1, \quad i\in\cM.
\end{cases}
\end{align*}
Since its generator is Lipschitz continuous, it admits a unique solution, denoted by $\big(P^{k, i}, \Lambda^{k, i}\big)_{i\in\cM}$.

It is easy to verify that
\begin{align}
(\underline P^i_t, \underline \Lambda^i_t)=\bigg(\Big(e^{-a'(T-t)}+\frac{1}{a'}(1-e^{-a'(T-t)})\Big)^{1-\gamma}, 0\bigg), \quad i\in\cM
\end{align}
is a solution to the following BSDE system:
\begin{align*}
\begin{cases}
d\underline P^i=-\Big[ -a'(1-\gamma) \underline P^i g(\underline P^{i})+(1-\gamma)(\underline P^i)^{-\frac{\gamma}{1-\gamma}}g(\underline P^{i})+\sum_{j=1}^\ell q^{ij}\underline P^j\Big]dt+(\underline\Lambda^i)'dW, \\
\underline P^i_T=1, \quad i\in\cM.
\end{cases}
\end{align*}
Since $\gamma<0$, we have the following estimates for $f^{k, i}(t, P, \Lambda)$ with $P>0$:
\begin{align*}
f^{k, i}(t, P, \Lambda)&\geq f^i(t, P, \Lambda)g(P)\\
&=\gamma \esssup_{(\pi, c)\in\Theta}\Big[-\frac{1-\gamma}{2}|\pi'\sigma^i|^2 P+\pi'(Pb^i+\sigma^i\Lambda)+\frac{c^\gamma}{\gamma}-Pc\Big]g(P)\\
&\geq \gamma \esssup_{(\pi, c)\in\R^m\times\R_+}\Big[-\frac{1-\gamma}{2}|\pi'\sigma^i|^2 P+\pi'(Pb^i+\sigma^i\Lambda)+\frac{c^\gamma}{\gamma}-Pc\Big]g(P)\\
&=\frac{\gamma}{2(1-\gamma)}\frac{1}{P}(Pb+\sigma\Lambda)'(\sigma\sigma')^{-1}(Pb+\sigma\Lambda)g(P)
+(1-\gamma)P^{-\frac{\gamma}{1-\gamma}}g(P).
\end{align*}
Hence, we deduce
\begin{align*}
&\quad f^{k, i}(t, \underline P^i, \underline \Lambda^i)-(\rho^i-\gamma r^i)\underline P^{i}g(\underline P^{i})+\sum_{j=1}^\ell q^{ij}\underline P^{j}\\
&\geq\frac{\gamma}{2(1-\gamma)}(b^i)'(\sigma\sigma')^{-1}b^i \underline P^i g(\underline P^{i})
+(1-\gamma)\underline P^{-\frac{\gamma}{1-\gamma}}g(\underline P^{i})-(\rho^i-\gamma r^i)\underline P^{i}g(\underline P^{i})+\sum_{j=1}^\ell q^{ij}\underline P^{j}\\
&\geq-a'(1-\gamma) \underline P^i g(\underline P^{i})
+(1-\gamma)\underline P^{-\frac{\gamma}{1-\gamma}}g(\underline P^{i})+\sum_{j=1}^\ell q^{ij}\underline P^{j}.
\end{align*}
By Lemma \ref{comparison}, we have
\[
a_1'= e^{-a'(1-\gamma)T}\leq\underline P^i_t\leq P^{k, i}_t, \quad i\in\cM.
\]

To give an upper bound, we notice that
\begin{align}
(\bar P^i, \bar \Lambda^i)=\Big(e^{a'(T-t)}+\frac{\varepsilon^\gamma}{a'}(e^{a'(T-t)}-1), 0\Big), \quad i\in\cM
\end{align}
is the unique solution to the following BSDE system:
\begin{align*}
\begin{cases}
d\bar P^i=-\Big[a'\bar P^i +\varepsilon^\gamma+\sum_{j=1}^\ell q^{ij}\bar P^j\Big]dt+(\bar\Lambda^i)'dW, \\
\bar P^i_T=1, \quad i\in\cM.
\end{cases}
\end{align*}
Also, $(0\mathbf{1}_m, \varepsilon)\in\Theta$ under Assumption \ref{assu0}, therefore for $k\geq-\gamma\varepsilon$, we have
\begin{align*}
f^{k, i}(t, P, \Lambda)&\leq \sup_{\tilde P\in\mathbb R, \tilde\Lambda\in\mathbb R^n}\Big[(\varepsilon^\gamma-\gamma \tilde P\varepsilon)g(\tilde P)-k|P-\tilde P|-k|\Lambda-\tilde\Lambda|\Big]\\
&\leq \sup_{\tilde P\in\mathbb R}\Big[(\varepsilon^\gamma-\gamma \tilde P\varepsilon)g(\tilde P)+\gamma\varepsilon|P-\tilde P|\Big]\\
&\leq \sup_{\tilde P\in\mathbb R}\Big[(\varepsilon^\gamma-\gamma \tilde P\varepsilon)g(\tilde P)+\gamma\varepsilon(|\tilde P|-|P|)\Big]\\
&\leq \varepsilon^\gamma-\gamma \varepsilon |P|,
\end{align*}
where the third and fourth inequalities are due to that $\gamma<0$ and $|xg(x)|\leq |x|$ for all $x\in\R$.
Hence, using \eqref{aprime2} and $\bar P^i>0$, we get
\begin{align*}
 f^{k, i}(t, \bar P^i, \bar \Lambda^i)-(\rho^i-\gamma r^i)\bar P^{i}g(\bar P^{i})+\sum_{j=1}^\ell q^{ij}\bar P^{j}
&\leq \varepsilon^\gamma-\gamma \varepsilon |\bar P^i|-(\rho^i-\gamma r^i)\bar P^{i}g(\bar P^{i})+\sum_{j=1}^\ell q^{ij}\bar P^{j}\\
&\leq a'\bar P^i +\varepsilon^\gamma+\sum_{j=1}^\ell q^{ij}\bar P^j.
\end{align*}
By Lemma \ref{comparison} again, we obtain the upper bound
\[
P^{k, i}_t\leq \bar P^i_t\leq e^{a'T}+\frac{\varepsilon^\gamma}{a'}(e^{a'T}-1)=a_2', \ t\in[0, T], \quad i\in\cM.
\]
\eof

Based on the above two theorems, we can provide the complete answer to problem \eqref{optm}.
\begin{theorem}
Suppose $\gamma\in(-\infty, 0)\cup(0, 1)$ and Assumptions \ref{assu1} and \ref{assu0} hold. Let
$(P^i, \Lambda^i)_{i\in\cM}$ be the unique solution to \eqref{P1}. Then the value function of the optimization problem \eqref{optm} is given by
$$V(x, i_{0})=\frac{1}{\gamma}x^\gamma P^{i_0}_0, $$
and $(\hat\pi, \hat c)\in\mathcal{U}$ is an optimal investment-consumption pair if and only if
\[
(\hat\pi_t, \hat c_t)\in\argmax_{(\pi, c)\in\Theta}\Big[-\frac{1-\gamma}{2}|\pi'\sigma_t^{\alpha_t}|^2 P_t^{\alpha_t}+\pi'(P_t^{\alpha_t}b_t^{\alpha_t}+\sigma_t^{\alpha_t}\Lambda_t^{\alpha_t})+\frac{c^\gamma}{\gamma}-P_t^{\alpha_t}c\Big].
\]
\end{theorem}

\begin{remark}
Let $G: \Omega\times[0,T]\times\mathbb{R}^m\times\mathbb{R}_+\rightarrow\mathbb{R}$ be a $\mathcal{P}\otimes\mathcal{B}(\mathbb{R}^m)\otimes\mathcal{B}(\mathbb{R}_+)$-measurable map,
where $\mathcal{P}$ is the $\sigma$-field of predictable sets of $\Omega\times[0,T]$, and $\mathcal{B}(\mathbb{R}^m)$ ($\mathcal{B}(\mathbb{R}_+)$) is the Borelian $\sigma$-algebra on $\mathbb{R}^m$ ($\mathbb{R}_+$).
By
\[
\argmax_{(\pi, c)\in\Theta}\Big[G(\omega, t, \pi, c)\Big],
\]
we denote the set of all predictable processes valued in $\Theta$ which attain the essential supremum of $G$ w.r.t $(\pi, c)\in\Theta$. From a measurable selection theorem (see e.g. \cite[Lemma 11]{HIM}, \cite[Corollary 18.14]{AB} or \cite[Proposition 2.4]{Sc}), the set is not empty. If $\Theta$ is further convex, then $\argmax_{(\pi, c)\in\Theta}\Big[G(t, \omega, \pi, c)\Big]$ degenerates to a singleton set.
\end{remark}

\begin{remark}
If $\Theta=\mathbb{R}^n\times\mathbb{R}_{+}$, then
\[
\hat\pi=\frac{1}{1-\gamma}(\sigma\sigma')^{-1}(b+\sigma \frac{\Lambda}{P}), \quad \hat c=P^{\frac{1}{\gamma-1}}.
\]
\end{remark}
\pf
For any $(\pi, c)\in\mathcal{U}$, applying It\^{o}'s formula to $\frac{1}{\gamma}e^{-\int_0^t\rho_s^{\alpha_s}ds}X_t^{\gamma}P_t^{\alpha_t}$, (here we use It\^{o}'s formula for Markovian chain; please refer to \cite[Lemma 4.3]{HLT}), we have
\begin{align*}
&\quad\frac{1}{\gamma}e^{-\int_0^t\rho_s^{\alpha_s}ds}X_t^{\gamma}P_t^{\alpha_t}
+\int_0^t \frac{1}{\gamma}e^{-\int_0^s\rho_u^{\alpha_u}du} X_s^{\gamma}c_s^{\gamma}ds\\
&=\frac{1}{\gamma}x^{\gamma}P_0^{i_0}+\int_0^t e^{-\int_0^s\rho_u^{\alpha_u}du} X_s^{\gamma}
\Big[\frac{1}{\gamma}c^\gamma -cP-\frac{1-\gamma}{2} P|\pi'\sigma|^2+\pi'(Pb+\sigma \Lambda)-\frac{f}{\gamma}\Big]ds\\
& \quad\;
+\int_0^t \frac{1}{\gamma}e^{-\int_0^s\rho_u^{\alpha_u}du} X_s^{\gamma} (\Lambda'+\gamma P\pi'\sigma)dW\\
& \quad\; +\int_0^t \frac{1}{\gamma} e^{-\int_0^s\rho_u^{\alpha_u}du} X_s^{\gamma}\sum_{j, j'\in\cM}(P^j-P^{j'})I_{\{\alpha_{s-}=j'\}} d\tilde N^{j'j}_s,
\end{align*}
where $(N^{j'j})_{j'j\in\mathcal{M}}$ are independent Poisson processes each with intensity $q^{j'j}$, and $\tilde N_t^{j'j}=N_t^{j'j}-q^{j'j}t, \ t\geq0$ are the corresponding compensated Poisson martingales under the filtration $\mathcal{F}$. And we drop the superscript $\alpha_s$ in the above integral.


Noting that for any $(\pi, c)\in\mathcal{U}$, the wealth process $X$ is continuous and strictly positive, hence bounded on $[0, T]$ away from $0$, a.s. Therefore, the stochastic integrals in the last equation are local martingales. Hence, there exists an increasing sequence of stopping times $(\tau_n)_{n\in\mathbb{N}}$ satisfying $\lim_{n\rightarrow\infty}\tau_n=T$, a.s. such that
\begin{align*}
&\quad\E\Big[\frac{1}{\gamma}e^{-\int_0^{\tau_n}\rho_s^{\alpha_s}ds}X_{\tau_n}^{\gamma}P_{\tau_n}^{\alpha_{\tau_n}}
+\int_0^{\tau_n} \frac{1}{\gamma}e^{-\int_0^s\rho_u^{\alpha_u}du} X_s^{\gamma}c_s^{\gamma}ds\Big]\\
&=\frac{1}{\gamma}x^{\gamma}P_0^{i_0}+\E\int_0^{\tau_n} e^{-\int_0^s\rho_u^{\alpha_u}du} X_s^{\gamma}
\Big[-\frac{1-\gamma}{2}P|\pi'\sigma|^2+\pi'(Pb+\sigma \Lambda)+\frac{1}{\gamma}c^\gamma -cP-\frac{f}{\gamma}\Big]ds\\
&\leq\frac{1}{\gamma}x^{\gamma}P_0^{i_0},
\end{align*}
thanks to the definition of $f$.

\textbf{Case $\gamma\in(0, 1)$}:
In this case, $\frac{1}{\gamma}e^{-\int_0^{\tau_n}\rho_s^{\alpha_s}ds}X_{\tau_n}^{\gamma}P_{\tau_n}^{\alpha_{\tau_n}}
+\int_0^{\tau_n} \frac{1}{\gamma}e^{-\int_0^s\rho_u^{\alpha_u}du} X_s^{\gamma}c_s^{\gamma}ds$ is bounded from below by $0$. Passing to limit and applying Fatou's lemma in the above inequality yields
\begin{align}\label{supermapower}
\E\Big[\frac{1}{\gamma}e^{-\int_0^T\rho_s^{\alpha_s}ds}X_T^{\gamma}
+\int_0^T \frac{1}{\gamma}e^{-\int_0^s\rho_u^{\alpha_u}du} X_s^{\gamma}c_s^{\gamma}ds\Big]
&\leq\frac{1}{\gamma}x^{\gamma}P_0^{i_0},
\end{align}
by virtue of the terminal condition
$P_T^{\alpha_T}=1$. And the inequality \eqref{supermapower} becomes an equality for some $(\hat\pi, \hat c)$ if and only if
\[
(\hat\pi, \hat c)\in \argmax_{(\pi, c)\in\Theta}\Big[-\frac{1-\gamma}{2}|\pi'\sigma^{\alpha_t}|^2 P+\pi'(P^{\alpha_t}b^{\alpha_t}+\sigma^{\alpha_t}\Lambda^{\alpha_t})
+\frac{c^\gamma}{\gamma}-P^{\alpha_t}c\Big].
\]
It remains to show that $(\hat\pi, \hat c)\in\mathcal{U}$. Recall that the unique solution of \eqref{P1} $(P^i, \Lambda^i)\in L^\infty_{\mathcal{F}^W}(0, T; \mathbb{R})\times \BMO$ for all $i\in\cM$.
From \eqref{Fboundeddomain1}, $|\hat\pi|\leq K(1+|Z^{\alpha_t}|)$, hence
$\hat\pi\in L^{2}_{\mathcal F}(0, T;\mathbb{R}^m)$.
Because $P$ is uniformly positive and bounded, we have
\begin{align*}
&\quad -\frac{1-\gamma}{2}|\pi'\sigma^i|^2 P+\pi'(Pb^i+\sigma^i\Lambda)+\frac{c^\gamma}{\gamma}-Pc\\
&\leq \esssup_{\pi\in\R^m}\Big[-\frac{1-\gamma}{2}|\pi'\sigma^i|^2 P+\pi'(Pb^i+\sigma^i\Lambda)\Big]+\frac{c^\gamma}{\gamma}-Pc\\
&=\frac{1}{2}\frac{\gamma}{1-\gamma}\frac{1}{P}(Pb+\sigma\Lambda)'(\sigma\sigma')^{-1}(Pb+\sigma\Lambda)
+\frac{c^\gamma}{\gamma}-Pc\\
&\leq K_{1}(1+|\Lambda|^2)+\frac{c^\gamma}{\gamma}-K_{1}c\\
&< 0,
\end{align*}
if $c> K(1+|\Lambda|^2)$ for sufficient large $K>0$. Hence $\hat c\leq K(1+|\Lambda|^2)$ so that $\int_0^T \hat cdt<\infty$, a.s. This implies $(\hat\pi, \hat c)\in\mathcal{U}$.

\textbf{Case $\gamma<0$}: For any admissible $(\pi, c)\in\mathcal{U}$, $(X^\gamma)_{0\leq t\leq T}$ belongs to class (D)
 and $\E\big[\int_0^T (cX)^\gamma dt\big]<\infty$. So by dominated convergence theorem, we also have \eqref{supermapower}.
 By a similar argument for the previous case $\gamma\in(0, 1)$, we can prove $\hat\pi\in L^{2}_{\mathcal F}(0, T;\mathbb{R}^m)$ and $\int_0^T \hat cdt<\infty$, a.s. It remains to prove $(\hat X^{\gamma}_t)_{0\leq t\leq T} \ \mbox{belongs to class (D) and } \E\big[\int_0^T (\hat c\hat X)^{\gamma}dt\big]<\infty$.

Since $\rho$ and $P$ are bounded, and
\begin{align*}
&\quad\E\bigg[e^{-\int_0^{\tau_n}\rho_s^{\alpha_s}ds}\hat X_{\tau_n}^{\gamma}P_{\tau_n}^{\alpha_{\tau_n}}
+\int_0^{\tau_n} e^{-\int_0^s\rho_u^{\alpha_u}du}\hat X_s^{\gamma}\hat c_s^{\gamma}ds\bigg]=x^{\gamma}P_0^{i_0},
\end{align*}
 Fatou's Lemma yields
\begin{align}\label{XTint:power}
\E\bigg[\hat X_T^{\gamma}+\int_0^T (\hat c\hat X)^\gamma dt\bigg]<\infty.
\end{align}

Denote $\theta^{\alpha_t}:=(\sigma^{\alpha_t})'(\sigma^{\alpha_t}(\sigma^{\alpha_t})')^{-1}$ for all $i\in\cM$.
By Girsanov's theorem,
\[
W^\mathbb{Q}_t:=W_t+\int_0^t\theta^{\alpha_s}_sds
\]
is a Brownian motion under the measure $\mathbb{Q}$ defined by
\[
\frac{d\mathbb{Q}}{d\mathbb{P}}=\mathcal{E}\Big(-\int_0^T\theta^{\alpha_s}_s dW_s\Big).
\]
Under the strategy $(\hat\pi, \hat c)$, we obtain from \eqref{wealth} that
\begin{align*}
\hat X_t=xe^{\int_0^t(r-\hat c)ds}\mathcal{E}\Big(\int_0^t\hat\pi'\sigma dW^{\mathbb{Q}}\Big).
\end{align*}
Denote
\[
J_t:=e^{-\int_0^t\hat cds}\mathcal{E}\Big(\int_0^t\hat\pi'\sigma dW^{\mathbb{Q}}\Big).
\]
Thus it is sufficient to prove $(J^{\gamma}_t)_{0\leq t\leq T}$ belongs to class (D), as $r$ is bounded.

It can be seen from \eqref{XTint:power} that
$\E[J_T^\gamma]<\infty$.
One obtains for every stopping time $\tau\leq T$,
\begin{align}\label{J1}
J_{\tau}^{\frac{\gamma}{2}}
&=\Big\{\E^{\mathbb{Q}}\Big[e^{-\int_0^{\tau}\hat cds}\mathcal{E}\Big(\int_0^{T}\hat\pi'\sigma dW^{\mathbb{Q}}\Big)\;\Big|\;\mathcal{F}_{\tau}\Big]\Big\}^{\frac{\gamma}{2}}\nonumber\\
&\leq\Big\{\E^{\mathbb{Q}}\Big[e^{-\int_0^{T}\hat cds}\mathcal{E}\Big(\int_0^{T}\hat\pi'\sigma dW^{\mathbb{Q}}\Big)\;\Big|\;\mathcal{F}_{\tau}\Big]\Big\}^{\frac{\gamma}{2}} =\Big\{\E^{\mathbb{Q}}\Big[J_T\;\Big|\;\mathcal{F}_{\tau}\Big]\Big\}^{\frac{\gamma}{2}} \leq\E^{\mathbb{Q}}\Big[J_T^{\frac{\gamma}{2}}\;\Big|\;\mathcal{F}_{\tau}\Big],
\end{align}
where we used that $\hat c> 0$ and $\gamma<0$ in the first inequality, and Jensen's inequality in the second one.
By the Cauchy-Schwarz inequality and the boundedness of $\theta$,
\begin{align}
\label{J2}
\E^{\mathbb{Q}}\Big[J_T^{\frac{\gamma}{2}}\;\Big|\;\mathcal{F}_{\tau}\Big]
&=\E\Big[J_T^{\frac{\gamma}{2}}\mathcal{E}\Big(-\int_{\tau}^T\theta_s^{\alpha_s}dW_s\Big)\;\Big|\;\mathcal{F}_{\tau}\Big]\nonumber\\
&\leq \Big\{\E\Big[J_T^{\gamma}\;\Big|\;\mathcal{F}_{\tau}\Big]\Big\}^{\frac{1}{2}}
\Big\{\E\Big[\mathcal{E}\Big(-\int_{\tau}^T\theta_s^{\alpha_s}dW_s\Big)^2\;\Big|\;\mathcal{F}_{\tau}\Big]\Big\}^{\frac{1}{2}}
\leq K \Big\{\E\Big[J_T^{\gamma}\;\Big|\;\mathcal{F}_{\tau}\Big]\Big\}^{\frac{1}{2}}.
\end{align}
Then \eqref{J1} and \eqref{J2} imply that $(J^{\gamma}_t)_{0\leq t\leq T}$ belongs to class (D).
\eof

\section{Logarithmic utility}\label{logU}
Let us now turn to the logarithmic utility $U(x)=\ln x, \ x>0$.
The admissible set is now defined as
\begin{align*}
\mathcal U&=\Big\{(\pi, c)\;\Big|\; \pi\in L^2_{\mathcal F}(0, T;\mathbb R^m), \ c\in L^1_{\mathcal F}(0, T;\mathbb{R}_{+}), \ (\pi_t(\omega), c_t(\omega))\in\Theta, \ a.e. , \ a.s, \\
&\qquad\qquad\quad (\ln X_t)_{0\leq t\leq T} \ \mbox{and} \ \Big(\int_0^t\ln(c_sX_s)ds\Big)_{0\leq t\leq T}\ \mbox{belong to class (D)} \Big\}.
\end{align*}
Clearly $\mathcal U$ would be empty if $c$ could not take positive value. Hence, it is reasonable to put the following assumption.
\begin{assumption}\label{assulog}
There exists $\varepsilon>0$ such that $(0\mathbf{1}_{m}, \varepsilon)\in\Theta$.
\end{assumption}

To solve problem \eqref{optm}, we introduce the following BSDE system:
\begin{align}\label{hlog}
\begin{cases}
dh_t^i=-(1-\rho^i h_t^i+\sum\limits_{j=1}^\ell q^{ij} h^j)dt+(\eta^i)'dW, \\
h_T^i=1, \quad i\in\cM.
\end{cases}
\end{align}
This is a linear BSDE system with bounded coefficients, so it has a unique solution $(h^i, \eta^i)_{i\in\cM}$ such that $(h^i, \eta^i)\in L^\infty_{\mathcal{F}^W}(0, T; \mathbb{R})\times \BMO$.

We claim that $h^i$ is uniformly positive. Indeed, under Assumption \ref{assu1}, there exists a constant $k>0$ such that $-\rho^i\geq -k$ for all $i\in\cM$. Notice that
\[
\Big(e^{-k(T-t)}+\frac{1-e^{-k(T-t)}}{k}, 0\Big)
\]
is a solution to the following BSDE:
\begin{align*}
\begin{cases}
dh_t=-(1-k h_t)dt+\eta'_{t}dW, \\
h_T=1,
\end{cases}
\end{align*}
From Lemma \ref{comparison}, we have $h^i_t\geq e^{-k(T-t)}+\frac{1-e^{-k(T-t)}}{k}\geq e^{-kT}$.

For the unique solution $(h^i, \eta^i)_{i\in\cM}$ of \eqref{hlog},
 consider the following BSDE system:
\begin{align}
\label{Plog}
\begin{cases}
dP^i=-\Big[f^i-\rho^iP^i+r^ih^i+\sum\limits_{j=1}^\ell q^{ij}P^j\Big]dt+(\Lambda^i)'dW \\
P^i_T=0, \quad i\in\cM,
\end{cases}
\end{align}
where $$f^i:=\esssup_{(\pi, c)\in\Theta}\Big[-\frac{1}{2}h^i|\pi'\sigma^i|^2+\pi'(h^ib^i+\sigma^i\eta^i)+\ln c-h^ic\Big].$$
As $(0\mathbf{1}_m, \varepsilon)\in\Theta$, and $h^i\in L^\infty_{\mathcal{F}^W}(0, T; \mathbb{R})$, we have
\begin{align}\label{floglower}
f^i\geq \ln\varepsilon-h^i\varepsilon.
\end{align}
On the other hand,
\begin{align*}
f^i
&\leq \esssup\limits_{(\pi, c)\in\mathbb{R}^m\times\mathbb{R}_+}\Big[-\frac{1}{2}h^i|\pi'\sigma^i|^2+\pi'(h^ib^i+\sigma^i\eta^i)
+\ln c-h^ic\Big]\\
&=\frac{1}{2h^i}(h^ib^i+\sigma^i\eta^i)'(\sigma^i(\sigma^i)')^{-1}(h^ib^i+\sigma^i\eta^i)
-1-\ln h^i.
\end{align*}
Thus
\begin{align*}
|f^i|\leq K(1+|\eta^i|^2).
\end{align*}
Whence \eqref{Plog} is an $\ell$-dimensional linear BSDE system with unbounded coefficients.
By \cite[Theorem 3.6]{HSX2}, the system of \eqref{Plog} admits a unique solution $(P^i, \Lambda^i)_{i\in\cM}$ such that
\[
(P^i, \Lambda^i)\in L^{\infty}_{\mathcal{F}^W}(0, T;\mathbb{R})\times\BMO, \ \mbox{for all} \quad i\in\cM.
\]

%
%
%

\begin{remark}
Set $(Y^i, Z^i)=(-\frac{P^i}{h^i}, \frac{P^i}{(h^i)^2}\eta^i-\frac{1}{h^i}\Lambda^i)$, then $(Y^i, Z^i)$ satisfies
\begin{align}
\label{Ylog}
\begin{cases}
dY^i=-\Big\{\essinf\limits_{(\pi, c)\in\Theta}\Big[\frac{1}{2}|\pi'\sigma^i|^2-\pi'(b^i+\frac{1}{h^i}\sigma^i\eta^i)
-\frac{\ln c}{h^i}+c\Big]-\frac{1}{h^i}(1+\sum\limits_{j=1}^\ell q^{ij}h^j)Y^i\\
\qquad\qquad\;+\frac{1}{h^i}(\red{\eta^i})'Z^i-r^i+\frac{1}{h^i}\sum\limits_{j=1}^\ell q^{ij}h^jY^j\Big\}dt+(Z^i)'dW, \\
Y^i_T=0, \quad i\in\cM.
\end{cases}
\end{align}
If there is no regime switching and the $\rho$ is a constant, then $\eta^i\equiv 0$ and \eqref{Ylog} degenerates to (4.6) in \cite{CH}.
\end{remark}

The following lemma can be found in Page 26 of \cite{Ka}.
\begin{lemma}
\label{exponpower}
Suppose $\phi\in \BMO$. 
Then
\begin{align*}
\E\left[\Big(\int_0^T|\phi_s|^2ds\Big)^j\right]<\infty, \quad \mbox{for any $j>0$}.
\end{align*}

\end{lemma}

\begin{theorem}
Suppose Assumptions \ref{assu1} and \ref{assulog} hold. Let $(h^i, \eta^i)_{i\in\cM}$ and
$(Y^i, Z^i)_{i\in\cM}$ be the unique solutions to \eqref{hlog} and \eqref{Ylog}, respectively. Then the value function of the optimization problem \eqref{optm} is given by
$$V(x, i_{0})=h_0^{i_0}\ln x+P_0^{i_0}, $$
and $(\hat\pi, \hat c)\in\mathcal{U}$ is an optimal consumption-investment pair if and only if
\[
(\hat\pi_t, \hat c_t)\in\argmax_{(\pi, c)\in\Theta}\Big[-\frac{1}{2}h_t^{\alpha_t}|\pi'\sigma_t^{\alpha_t}|^2
+\pi'(h_t^{\alpha_t}b_t^{\alpha_t}+\sigma_t^{\alpha_t}\eta_t^{\alpha_t})
+\ln c-h_t^{\alpha_t}c\Big].
\]
\end{theorem}
\begin{remark}
If $\Theta=\mathbb{R}^m\times\mathbb{R}_+$, then the optimal investment-consumption pair is uniquely given by
$$\hat\pi_t=(\sigma^{\alpha_t}(\sigma^{\alpha_t})')^{-1}(b^{\alpha_t}+\frac{1}{h^{\alpha_t}}\sigma^{\alpha_t}\eta^{\alpha_t}), \quad \hat c_t=\frac{1}{h_t^{\alpha_t}}.$$
\end{remark}
\pf
For any $(\pi, c)\in\mathcal{U}$, applying It\^{o}'s formula to $e^{-\int_0^t\rho_s^{\alpha_s}ds}(h_t^{\alpha_t}\ln X_t+P_t^{\alpha_t})$, we have
\begin{align*}
&\quad e^{-\int_0^t\rho_s^{\alpha_s}ds}(h_t^{\alpha_t}\ln X_t+P_t^{\alpha_t})+\int_0^t e^{-\int_0^s\rho_u^{\alpha_u}du}\ln(c_sX_s)ds\\
&=h_0^{i_0}\ln x +P_0^{i_0}+\int_0^t e^{-\int_0^s\rho_u^{\alpha_u}du}\Big[\ln c-ch-\frac{1}{2}|\pi'\sigma|^2+\pi'(h\eta+\sigma\eta)-f^{\alpha_s}\Big]ds\\
&\qquad+\int_0^t e^{-\int_0^t\rho_s^{\alpha_s}ds}\Big[h\sigma'\pi+\eta\ln X+\Lambda\Big]'dW\\
&\qquad+\int_0^t e^{-\int_0^t\rho_s^{\alpha_s}ds}\sum_{j, j'\in\cM}[(h^j\ln X+P^j)-(h^{j'}\ln X+P^{j'})]I_{\{\alpha_{s-}=j'\}} d\tilde N^{j'j}_s.
\end{align*}
Then by the definition of $f^i$,
the process
\[
e^{-\int_0^t\rho_s^{\alpha_s}ds}(h_t^{\alpha_t}\ln X_t+P_t^{\alpha_t})+\int_0^t e^{-\int_0^s\rho_u^{\alpha_u}du}\ln(c_sX_s)ds
\]
is a local supermartingale.
Hence, there exists a sequence of increasing stopping times $(\tau_n)_{n\in\mathbb{N}}$ satisfying $\lim_{n\rightarrow\infty}\tau_n=T$, a.s. such that
\begin{align*}
\E\Big[e^{-\int_0^{\tau_n}\rho_s^{\alpha_s}ds}(h_{\tau_n}^{\alpha_{\tau_n}}\ln X_{\tau_n}+P_{\tau_n}^{\alpha_{\tau_n}})+\int_0^{\tau_n} e^{-\int_0^s\rho_u^{\alpha_u}du}\ln(c_sX_s)ds\Big]\leq h_0^{i_0}\ln x+P_0^{i_0}.
\end{align*}
Recalling the admissibility of $(\pi, c)$, we deduce from the dominated convergence theorem that
\begin{align}\label{logclassd}
\E\Big[e^{-\int_0^{T}\rho_s^{\alpha_s}ds}(h_{T}^{\alpha_{T}}\ln X_{T}+P_{T}^{\alpha_{T}})+\int_0^{T} e^{-\int_0^s\rho_u^{\alpha_u}du}\ln(c_sX_s)ds\Big]\leq h_0^{i_0}\ln x+P_0^{i_0},
\end{align}
with the equality holds for some $(\hat\pi, \hat c)$
if and only if
\[
(\hat\pi_t, \hat c_t)\in\argmax_{(\pi, c)\in\Theta}\Big[-\frac{1}{2}h^i|\pi'\sigma^i|^2+\pi'(h^ib^i+\sigma^i\eta^i)
+\ln c-h^ic\Big].
\]
It remains to prove $(\hat\pi, \hat c)\in\mathcal{U}$ for such a pair.
\footnote{
Compared with \cite{CH}, a new, non-zero term $\eta\in\BMO$ emerges (In \cite{CH}, $\rho$ is a constant, hence $\eta=0$). So we need conditions like class (D) to apply the dominated convergence theorem in \eqref{logclassd}. But this brings difficulty in proving the admissibility of $(\hat\pi, \hat c)$.}

From \eqref{floglower} and $h^i\in L^{\infty}_{\mathcal{F}^W}(0, T;\mathbb{R})$, $f^i\geq -K_1$ for some $K_1>0$.
On the other hand, Assumption \ref{assu1} and $h^i\geq e^{-kT}$ imply
\begin{align*}
&\quad -\frac{1}{2}h^i|\pi'\sigma^i|^2+\pi'(h^ib^i+\sigma^i\eta^i)+\ln c-h^ic\\
&\leq -\frac{1}{2}h^i|\pi'\sigma^i|^2+\pi'(h^ib^i+\sigma^i\eta^i)+\esssup_{c\in\R_+}(\ln c-h^ic)\\
&=-\frac{1}{2}h^i|\pi'\sigma^i|^2+\pi'(h^ib^i+\sigma^i\eta^i)-\ln h^i-1\\
&\leq -K_2|\pi|^2+K_3|\pi|(1+|\eta^i|)+K_4,
\end{align*}
for some positive constants $K_2, \ K_3, \ K_4$. Thus
\begin{align*}
&\quad -\frac{1}{2}h^i|\pi'\sigma^i|^2+\pi'(h^ib^i+\sigma^i\eta^i)+\ln c-h^ic <-K_1\leq f^i,
\end{align*}
if $|\pi|> K(1+|\eta^i|)$ with $K$ being sufficient large, which implies
\begin{align}\label{pilog}
|\hat\pi|\leq K(1+|\eta^i|).
\end{align}
Similarly, since $h^i$ is uniformly positive, we have for some $K_5, K_{6}>0$,
\begin{align*}
&\quad -\frac{1}{2}h^i|\pi'\sigma^i|^2+\pi'(h^ib^i+\sigma^i\eta^i)+\ln c-h^ic\\
&\leq \esssup_{\pi\in\R^m}\Big[-\frac{1}{2}h^i|\pi'\sigma^i|^2+\pi'(h^ib^i+\sigma^i\eta^i)\Big]+\ln c-K_{5} c\\
&=\frac{1}{2h^i}(h^ib^i+\sigma^i\eta^i)'(\sigma^i(\sigma^i)')^{-1}(h^ib^i+\sigma^i\eta^i)+\ln c-K_{5} c\\
&\leq K_{6}(1+|\eta^i|^2)+\ln c-K_{5}c\\
&<-K_1\\
&\leq f^i,
\end{align*}
if $c\geq K(1+|\eta^i|^2)$ or $\ln c\leq -K(1+|\eta^i|^2)$ for sufficient large $K$. Therefore,
\begin{align}\label{clog}
e^{-K(1+|\eta^i|^2)}\leq \hat c\leq K(1+|\eta^i|^2),
\end{align}
and consequently,
\begin{align}\label{clog2}
|\ln\hat c|\leq K(1+|\eta^i|^2).
\end{align}

From \eqref{wealth},
\begin{align*}
\ln\hat X_t&=\ln(x)+\int_0^t (r+\hat\pi'b-\hat c-\frac{1}{2}|\hat\pi'\sigma|^2)ds+\int_0^t \hat\pi'\sigma dW.
\end{align*}
It follows from the boundedness of coefficients and the AM-GM inequality that
\begin{align*}
\sup_{t\leq T}|\ln\hat X_t|&\leq K+K\int_0^T (|\hat c|+ |\hat\pi|^2)ds+K\sup_{t\leq T}\bigg|\int_0^t \hat\pi'\sigma dW\Big|.
\end{align*}
Using the AM-GM inequality, the Burkholder-Davis-Gundy inequality, \eqref{pilog}, \eqref{clog}, and Lemma \ref{exponpower}, we have
\begin{align*}
\E\Big[\sup_{\tau\leq T}(\ln\hat X_\tau)^2\Big]
& \leq K+K\E\bigg[\Big(\int_0^{T}(|\hat c|+|\hat\pi|^{2})ds\Big)^2+\sup_{t\leq T} \bigg|\int_0^t\hat\pi'\sigma dW\bigg|^2\bigg]\\
&\leq K+K\E\bigg[\Big(\int_0^{T}(1+|\eta^i|^2)ds\Big)^2+\int_0^T|\hat\pi'\sigma|^{2}ds \bigg]\\
&\leq K+K\E\bigg[\Big(\int_0^{T}(1+|\eta^i|^2)ds\Big)^2+\int_0^{T}(1+|\eta^i|^2)ds \bigg]\\
&\leq K.
\end{align*}
This shows that $(\ln X_t)_{0\leq t\leq T}$ belongs to class (D).
By the same argument and using \eqref{clog2}, we can prove that $\Big(\int_0^t\ln(c_sX_s)ds\Big)_{0\leq t\leq T}$ belongs to class (D) as well.
\eof

\section{Exponential utility}\label{expU}
In this case, we do not consider proportional consumption-investment strategies. Instead,
let $\pi_{j, t}$ be the \emph{amount} invest in the $j$th risky asset, $j=1, \ldots, m$, and let $ c_t$ be the \emph{amount} to consume.
Then the investor's self-financing wealth process $X$ corresponding to a consumption-investment strategy $(\pi,c)$ is the unique strong solution of the SDE:
\begin{align}
\label{wealthex}
\begin{cases}
dX_t=(r_tX_t+\pi_t'b_t^{\alpha_t}-c_t )dt+\pi_t'\sigma_t^{\alpha_t}dW_t, \\
X_0=x, \ \alpha_0=i_0.
\end{cases}
\end{align}

Let $\Pi$ be a given closed nonempty set in $\mathbb{R}^m$ with $0\in\Pi$ to represent the constraint set for portfolios. We do not put constraint on consumptions, so $-c_{t}$ can be interpreted as income if $c_{t}<0$. The admissible consumption-investment set is defined as
\begin{align*}
\mathcal U&=\bigg\{(\pi, c)\;\bigg|\; \int_0^T(|\pi_t|^2+|c_t|)dt<\infty, \ (\pi_t(\omega), c_t(\omega))\in\Pi\times\mathbb{R}, \ a.e. , \ a.s, \\
&\qquad\qquad\quad\; \mbox{and} \ (e^{-\beta h_tX_t^{\pi, c}})_{0\leq t\leq T} \ \mbox{ belongs to class (D)} \bigg\},
\end{align*}
where the process $h_t$ will be specified case by case in the following subsections.

\par
The investor's problem is to maximize
\begin{align}
J(x, i_0;\pi, c):&=\E\bigg[\int_0^T e^{-\int_0^t\rho_s^{\alpha_s}ds}U(c_t)dt+e^{-\int_0^T\rho_s^{\alpha_s}ds}U(X_T)\bigg]%
, \quad \mathrm{ s.t.} \;\;
(\pi, c)\in \mathcal{U},
\label{optmexp}%
\end{align}
and determine the value function $$V(x, i_0):=\sup_{(\pi, c)\in\mathcal{U}}J(x, i_0;\pi, c), $$
where the utility is exponential and given by \[U(x)=-\exp(-\beta x), \quad x\in\mathbb{R},\quad \beta>0.\]

It turns out that whether the interest rate is deterministic or random is crucial, we have to deal with these two cases separately.

\subsection{Case I: Deterministic interest rate}
 \begin{assumption} \label{assu1exp}
The process $r$ is a deterministic bounded measurable function of $t$. For all $i\in\cM$,
\begin{align*}
\mu^{i}\in L_{\mathbb{F}^W}^\infty(0, T;\mathbb{R}^m), \quad
\sigma^{i}\in L_{\mathbb{F}^W}^\infty(0, T;\mathbb{R}^{m\times n}), \quad
\rho^i\in L_{\mathbb{F}^W}^\infty(0, T;\mathbb{R}), \
\end{align*}
and $\sigma^{i}(\sigma^{i})'\geq\delta I_{m}$ with some constant $\delta>0$.
\end{assumption}

Define
\begin{align}\label{hexp}
h_t=\Big(e^{-\int_t^T r_sds}+\int_t^T e^{-\int_t^s r_udu}ds\Big)^{-1}>0.
\end{align}
Then under Assumption \ref{assu1exp}, $h$ is the unique bounded uniformly positive solution to the following ODE:
\begin{align*}
\begin{cases}
dh_t=-h_t(r_t-h_t)dt, \\
 h_T=1.
\end{cases}
\end{align*}
Consider the following BSDE system:
\begin{align}
\label{YexpP}
\begin{cases}
dP^i=-\Big\{f^i(P^i, \Lambda^i)-hP^i\ln P^i-\rho^iP^i+h(1-\ln h)P^i+\sum_{j=1}^\ell q^{ij}P^j\Big\}dt+(\Lambda^i)'dW, \\
P^i_T=1, \quad i\in\cM,
\end{cases}
\end{align}
where
\begin{align*}
f^i(P, \Lambda):=\beta h\essinf_{\pi\in\Pi}\Big[\frac{1}{2}\beta hP^i|\pi'\sigma^i|^2-\pi'(P^ib^i+\sigma^i\Lambda^i)\Big].
\end{align*}
Thanks to $0\in\Pi$, we have $f^{i}\leq 0$.

Solutions to \eqref{YexpP} are defined as in Definition \ref{def}.

\begin{theorem}\label{Th:Pexp:deter}
Under Assumption \ref{assu1exp}, there is a unique uniformly positive solution $(P^i, \Lambda^i)_{i\in\cM}$ to the BSDE \eqref{YexpP}.
\end{theorem}
\pf
We will only present how to find the uniformly bounds. Details are left to interested readers.

Let $a_1>0$, $a_2>0$, $a_3>0$ be three constants such that
\begin{align}\label{a3}
&-\rho^i+h(1-\ln h)\leq a_1, \nonumber\\
&\frac{1}{2}(b^i)'(\sigma^i(\sigma^i)')^{-1}b^i-\rho^i+h(1-\ln h)\geq -a_2, \nonumber\\
&h\geq a_3>0, \ \forall i\in\cM.
\end{align}
Denote $a:=e^{a_1 T}$ and $\epsilon:=e^{-\frac{a_2}{a_3}(1-e^{-a_3T})}$. Note $\epsilon<1<a$.

Let $g:\R \to[0, 1]$ be a smooth truncation function satisfying $g(x)=0$ for $x\in (-\infty, \epsilon/2]\cup[2a, \infty)$, and $g(x)=1$ for $x\in[\epsilon, a]$.

For $k\geq1$, $(t, P, \Lambda)\in[0, T]\times\mathbb R\times\mathbb{R}^n$, $ i\in\cM$, define
\[
f^{k, i}(t, P, \Lambda)=\sup_{\tilde P\in\mathbb R, \tilde\Lambda\in\mathbb R^n}\Big[f^i (t, \tilde P, \tilde\Lambda)g(\tilde P)-k|P-\tilde P|-k|\Lambda-\tilde\Lambda|\Big].
\]
Let $\big(P^{k, i}, \Lambda^{k, i}\big)_{i\in\cM}$ denote the unique solution to the following $\ell$-dimensional BSDE system with a Lipschitz generator
\begin{align*}
\begin{cases}
dP^{k, i}=-\Big[f^{k, i}(P^{k, i}, \Lambda^{k, i, +})-hP^{k, i}(\ln P^{k, i})g(P^{k, i})-\rho^iP^{k, i}g(P^{k, i})\\
 \qquad\qquad\quad+h(1-\ln h)P^{k, i}g(P^{k, i})+\sum_{j=1}^\ell q^{ij}P^{k, j}\Big]dt+(\Lambda^{k, i})^\top dW, \\
P^{k, i}_T=1, \quad i\in\cM.
\end{cases}
\end{align*}

It is direct to verify that
\begin{align}
(\underline P^i_t, \underline \Lambda^i_t)=\Big(e^{-\frac{a_2}{a_3}(1-e^{-a_3(T-t)})}, 0\Big), \quad i\in\cM
\end{align}
is the unique solution to the following BSDE system:
\begin{align*}
\begin{cases}
d\underline P^i=-\Big[-a_2\underline P^ig(\underline P^{i})-a_3 \underline P^i (\ln \underline P^i)g(\underline P^{i}) +\sum_{j=1}^\ell q^{ij}\underline P^j\Big]dt+(\underline\Lambda^i)'dW, \\
\underline P^i_T=1, \quad i\in\cM.
\end{cases}
\end{align*}
Notice that $\underline P^i \in(0, 1)$, and
\begin{align*}
f^{k, i}(t, P, \Lambda)&\geq f^i(t, P, \Lambda)g(P)\\
&= \beta h\essinf_{\pi\in\Pi}\Big[\frac{1}{2}\beta hP^i|\pi'\sigma^i|^2-\pi'(P^ib^i+\sigma^i\Lambda^i)\Big]g(P)\\
&\geq \beta h\essinf_{\pi\in\mathbb{R}^m}\Big[\frac{1}{2}\beta hP^i|\pi'\sigma^i|^2-\pi'(P^ib^i+\sigma^i\Lambda^i)\Big]g(P)\\
&= \frac{1}{2}(b^i)'(\sigma^i(\sigma^i)')^{-1}b^iPg(P)+ (b^i)'(\sigma^i(\sigma^i)')^{-1}\sigma\Lambda g(P)\\
&\qquad+ \frac{1}{2P}\Lambda'(\sigma^i)'(\sigma^i(\sigma^i)')^{-1}\sigma^i\Lambda g(P),
\end{align*}
we deduce from \eqref{a3} that
\begin{align*}
&\quad f^{k, i}(\underline P^i, \underline \Lambda^i)-h\underline P^i (\ln \underline P^i)g(\underline P^{i})-\rho^i\underline P^ig(\underline P^{i})+h(1-\ln h)\underline P^ig(\underline P^{i})+\sum_{j=1}^\ell q^{ij}\underline P^i\\
&\geq \frac{1}{2}(b^i)'(\sigma^i(\sigma^i)')^{-1}b^i\underline P^ig(\underline P^{i})-h\underline P^i(\ln \underline P^i) g(\underline P^{i})-\rho^i\underline P^ig(\underline P^{i})+h(1-\ln h)\underline P^ig(\underline P^{i})+\sum_{j=1}^\ell q^{ij}\underline P^i\\
&\geq-a_2\underline P^ig(\underline P^{i})-a_3 \underline P^i(\ln \underline P^i)g(\underline P^{i}) +\sum_{j=1}^\ell q^{ij}\underline P^j.
\end{align*}
By Lemma \ref{comparison}, we have
\begin{align}
P^{k, i}_t\geq\underline P^i_t\geq e^{-\frac{a_2}{a_3}(1-e^{-a_3T})}=\epsilon, \quad i\in\cM.
\end{align}
This gives a uniform lower bound.

To give a uniform upper bound, we notice
\begin{align}
(\bar P^i_{t}, \bar \Lambda^i_{t})=(e^{a_1(T-t)}, 0), \quad i\in\cM,
\end{align}
is a solution to the following BSDE system:
\begin{align*}
\begin{cases}
d \bar P^i=-\Big[a_1\bar P^ig(\bar P^{i})+\sum_{j=1}^\ell q^{ij}\bar P^j\Big]dt+(\bar\Lambda^i)'dW, \\
\bar P^i_T=1, \quad i\in\cM.
\end{cases}
\end{align*}
Since $0\in\Pi$, we have $f^i(t, P, \Lambda)\leq 0$, and hence $f^{k, i}(t, P, \Lambda)\leq 0$. Notice that $\bar P^i_t\geq 1$, thus
we have
\begin{align*}
&\quad f^{k, i}(\bar P^i, \bar \Lambda^i)-h\bar P^i (\ln \bar P^i)g(\bar P^{i})-\rho^i\bar P^ig(\bar P^{i})+h(1-\ln h)\bar P^ig(\bar P^{i})+\sum_{j=1}^\ell q^{ij}\bar P^j\\
&\leq a_1\bar P^ig(\bar P^{i})+\sum_{j=1}^\ell q^{ij}\bar P^j.
\end{align*}
By Lemma \ref{comparison}, we have
\begin{align}
P^{k, i}_t\leq\bar P^i_t\leq e^{a_1T}=a, \quad i\in\cM.
\end{align}
This gives a uniform upper bound.
\eof

Set $Y^i=-\frac{1}{\beta}\ln P^i$, $Z^i=-\frac{\Lambda^i}{\beta P^i}$, then $(Y^i, Z^i)$ is the unique solution to
\begin{align}
\label{Yexp}
\begin{cases}
dY^i=-\Big[F^i(Z^i)-hY^i-\frac{\beta}{2}|Z^i|^2+\frac{\rho^i}{\beta}-\frac{h}{\beta}(1-\ln h)\\
\qquad\qquad-\frac{1}{\beta}\sum_{j=1}^\ell q^{ij}\Big(e^{-\beta(Y^j-Y^i)}-1\Big)\Big]dt+(Z^i)'dW, \\
Y^i_T=0, \quad i\in\cM,
\end{cases}
\end{align}
where
\[
F^i(z):=h\esssup_{\pi\in\Pi}\Big[-\frac{1}{2}\beta h|\pi'\sigma^i|^2+\pi'(b^i-\beta\sigma^i z)\Big],\quad z\in\R^n.
\]

From above, we immediately have
\begin{theorem}\label{Th:Yexp:deter}
Under Assumption \ref{assu1exp}, there is a unique solution $(Y^i, Z^i)_{i\in\cM}$ to the BSDE \eqref{Yexp} such that $(Y^i, Z^i)\in L^\infty_{\mathcal{F}^W}(0, T; \mathbb{R})\times \BMO$ for all $i\in\cM$.
\end{theorem}

Now we can present the solution to problem \eqref{optmexp}.
\begin{theorem}
Suppose Assumption \ref{assu1exp} holds. Let $h$ be defined in \eqref{hexp} and
$(Y^i, Z^i)_{i\in\cM}$ be the unique solution to \eqref{Yexp}. Then the value function of the optimization problem \eqref{optmexp} is given by
\[V(x, i_{0})=-e^{-\beta(h_0 x+Y_0^{i_0})}, \]
and the optimal portfolio is
\begin{align*}
\hat\pi_t\in\argmax_{\pi\in\Pi}\Big[-\frac{1}{2}\beta h_t|\pi'\sigma_t^{\alpha_t}|^2+\pi'(b_t^{\alpha_t}-\beta\sigma_t^{\alpha_t} Z_{t}^{\alpha_t})\Big],
\end{align*}
the optimal consumption is
\begin{align*}
\hat c_t=h_t X_t+Y_t^{\alpha_t}-\frac{1}{\beta}\ln h_t.
\end{align*}
\end{theorem}
\pf
For any $(\pi, c)\in\mathcal{U}$, applying It\^{o}'s formula to $-e^{-\int_0^t\rho_s^{\alpha_s}ds}e^{-\beta(h_tX_t+Y_t^{\alpha_t})}$, we have
\begin{align*}
&\quad -e^{-\int_0^t\rho_s^{\alpha_s}ds}e^{-\beta(h_tX_t+Y_t^{\alpha_t})}
-\int_0^t e^{-\int_0^s\rho_s^{\alpha_u}du}e^{-\beta c_s}ds\\
&=-e^{-\beta(h_0x+Y_0^{i_0})}+
\int_0^t\beta e^{-\int_0^s\rho_u^{\alpha_u}du}e^{-\beta(h_sX_s+Y_s^{\alpha_s})}
\Big[-\frac{1}{\beta}e^{\beta (hX+Y)}e^{-\beta c}-hc\\
&\qquad\qquad\qquad-\Big(\frac{h}{\beta}(\ln h-1)-h^2X-hY\Big) +h\Big(-\frac{1}{2}\beta h|\pi'\sigma^i|^2+\pi'(b^i-\beta\sigma^i z)\Big)-F^i(Z)\Big]ds\\
&\qquad+\int_0^t \beta e^{-\int_0^s\rho_u^{\alpha_u}du}e^{-\beta(h_sX_s+Y_s^{\alpha_s})}(h\pi'\sigma+Z')dW\\
&\qquad+\int_0^t e^{-\int_0^s\rho_u^{\alpha_u}du}e^{-\beta h_sX_s} \sum_{j, j'\in\cM}(e^{-\beta Y^j}-e^{-\beta Y^i})I_{\{\alpha_{s-}=j'\}}d\tilde N^{j'j}_s,
\end{align*}
where $(N^{j'j})_{j'j\in\mathcal{M}}$ are independent Poisson processes each with intensity $q^{j'j}$, and $\tilde N_t^{j'j}=N_t^{j'j}-q^{j'j}t, \ t\geq0$ are the corresponding compensated Poisson martingales under the filtration $\mathcal{F}$. And we drop the superscript $\alpha_s$ in the above integral.

Noting that for any $(\pi, c)\in\mathcal{U}$, the wealth process $X$ is continuous, hence bounded on $[0, T]$ a.s. Therefore, the stochastic integrals in the last equation are local martingales. Hence, there exists a sequence of increasing stopping times $(\tau_n)_{n\in\mathbb{N}}$ satisfying $\lim_{n\rightarrow\infty}\tau_n=T$, a.s. such that
\begin{align*}
&\quad \E\Big[-e^{-\int_0^{\tau_n}\rho_s^{\alpha_s}ds}e^{-\beta(h_{\tau_n}X_{\tau_n}+Y_{\tau_n}^{\alpha_{\tau_n}})}
-\int_0^{\tau_n} e^{-\int_0^s\rho_s^{\alpha_u}du}e^{-\beta c_s}ds\Big]\\
&=-e^{-\beta(h_0x+Y_0^{i_0})}+
\int_0^{\tau_n}\beta e^{-\int_0^s\rho_u^{\alpha_u}du}e^{-\beta(h_sX_s+Y_s^{\alpha_s})}
\Big[-\frac{1}{\beta}e^{\beta (hX+Y)}e^{-\beta c}-hc\\
&\qquad\qquad\qquad-\Big(\frac{h}{\beta}(\ln h-1)-h^2X-hY\Big) +h\Big(-\frac{1}{2}\beta h|\pi'\sigma^i|^2+\pi'(b^i-\beta\sigma^i z)\Big)-F^i(Z)\Big]ds.
\end{align*}
Because $h,\beta>0$,
$$\esssup_{c\in\R }\Big[-\frac{1}{\beta}e^{\beta (hX+Y)}e^{-\beta c}-hc\Big]=\frac{h}{\beta}(\ln h-1)-h^2X-hY,$$
we obtain from the definition of $F^i$ that
\begin{align*} \E\bigg[-e^{-\int_0^{\tau_n}\rho_s^{\alpha_s}ds}e^{-\beta(h_{\tau_n}X_{\tau_n}+Y_{\tau_n}^{\alpha_{\tau_n}})}\bigg]
-\E\bigg[\int_0^{\tau_n} e^{-\int_0^s\rho_s^{\alpha_u}du}e^{-\beta c_s}ds\bigg]\leq -e^{-\beta(h_0x+Y_0^{i_0})},
\end{align*}
for any $(c, \pi)\in\mathcal{U}$. As $n\to\infty$, the second expectation is convergent by the monotone convergence theorem; while the first expectation, because $(e^{-\beta h_tX_t^{\pi, c}})_{0\leq t\leq T}$ belongs to class (D), $Y^{i}$, $i\in\cM$, are bounded, is also convergent by the dominated convergence theorem. Hence
\begin{align} \label{supermader}
\E\bigg[-e^{-\int_0^{T}\rho_s^{\alpha_s}ds}e^{-\beta X_{T}}
-\int_0^{T} e^{-\int_0^s\rho_s^{\alpha_u}du}e^{-\beta c_s}ds\bigg]
\leq -e^{-\beta(h_0x+Y_0^{i_0})},
\end{align}
by virtue of the terminal conditions $h_T=1$ and $Y_T=0$. Furthermore, the inequality \eqref{supermader} becomes an equality for some $(\hat\pi, \hat c)$ if and only if
\begin{align*}
\hat c_t&=\argmax_{c\in\R}\Big[-\frac{1}{\beta}e^{\beta (hX+Y^{\alpha})}e^{-\beta c}-hc\Big] =-\frac{1}{\beta}\ln h +hX+Y^{\alpha},
\end{align*}
and
\begin{align*}
\hat\pi_t\in\argmax_{\pi\in\Pi}\Big[-\frac{1}{2}\beta h|\pi'\sigma^{\alpha}|^2+\pi'(b^{\alpha}-\beta\sigma^{\alpha} Z^{\alpha})\Big].
\end{align*}

The admissibility issue $(\hat\pi, \hat c)\in\mathcal{U}$ is covered in a more general case; please refer to the proof of Theorem \ref{Th:opt}.
\eof


\subsection{Case II: Random interest rate}
\begin{assumption} \label{assu:r:ran}
For all $i\in\cM$,
\begin{align*}
\rho^i\in L_{\mathbb{F}^W}^\infty(0, T;\mathbb{R}).
\end{align*}
Furthermore, the processes $r, \ \mu, \ \sigma$ are independent of the regime $i\in\cM$ and such that
\begin{align*}
r\in L^{\infty}_{\mathcal{F}^W}(0, T;\mathbb{R}), \ \mu\in L_{\mathbb{F}^W}^\infty(0, T;\mathbb{R}^m), \
\sigma\in L_{\mathbb{F}^W}^\infty(0, T;\mathbb{R}^{m\times n}),
\end{align*}
and $\sigma\sigma'\geq\delta I_{m}$ with some constant $\delta>0$.
\end{assumption}
\begin{assumption}\label{assu:r:ran2}
$m=n$ and $\Pi=\mathbb{R}^m$, i.e. neither the portfolio nor the consumption has constraints.
\end{assumption}
Under Assumptions \ref{assu:r:ran} and \ref{assu:r:ran2}, $\sigma$ is a square invertible matrix.

Consider the following BSDE system:
\begin{align}\label{hexp:ran}
\begin{cases}
dh_t=-[rh-h^2-\eta'\sigma^{-1}b-\frac{1}{h}|\eta|^2]dt+\eta'dW, \\
 h_T=1.
\end{cases}
\end{align}
\begin{theorem}\label{Th:exp:h}
Under Assumptions \ref{assu:r:ran} and \ref{assu:r:ran2}, \eqref{hexp:ran} admits a unique solution $(h, \eta)\in L^\infty_{\mathcal{F}^W}(0, T; \mathbb{R})\times \BMO$ such that $h\geq \delta$ for some positive constant $\delta$.
\end{theorem}
\pf
Let $(p, q)\in L^\infty_{\mathcal{F}^W}(0, T; \mathbb{R})\times \BMO$ denote
unique solution to the following linear BSDE system:
\begin{align}\label{hexp:p}
\begin{cases}
dp_t=-[1-rp-q'\sigma^{-1}b]dt+q'dW, \\
 p_T=1.
\end{cases}
\end{align}
Indeed, $p$ admits the explicit representation:
\begin{align*}
p_t=\widetilde{\E}\bigg[e^{-\int_t^Tr_sds}+\int_t^T e^{-\int_t^sr_udu}ds \;\bigg|\;\mathcal{F}^W_t\bigg],
\end{align*}
where $\widetilde{\E}$ is the expectation with respect to the probability measure $\widetilde {\mathbb{P}}$ defined by
\[
\frac{d\widetilde {\mathbb{P}}}{d\mathbb{P}}=\mathcal{E}\bigg(\int_0^T(\sigma^{-1}b)'dW\bigg).
\]
It follows from Assumption \ref{assu:r:ran} that $\delta_1\leq p\leq \delta_2$ for some $0<\delta_1<\delta_2$. Therefore $(h, \eta):=(\frac{1}{p}, -\frac{q}{p^2})$ is well defined and $(h, \eta)\in L^\infty_{\mathcal{F}^W}(0, T; \mathbb{R})\times \BMO$ with $h\geq \frac{1}{\delta_2}>0$. It can be directly verified, using It\^{o}'s formula, that $(h, \eta)$ is a solution of \eqref{hexp:ran}. As for uniqueness, observe that if $(h, \eta)$ is a solution of \eqref{hexp:ran} with $h\geq \delta$ for some $\delta>0$, then $(p, q^i)=(\frac{1}{h}, -\frac{\eta}{h^2})$ is a solution to \eqref{hexp:p}. Uniqueness follows from the fact that \eqref{hexp:p} has a unique solution.
\eof

\bigskip

Let $(h, \eta)$ denote the unique solution to \eqref{hexp:ran}. Consider the following BSDE system:
\begin{align}
\label{Yexp:ran}
\begin{cases}
dY^i=-\Big[-hY^i-(\sigma^{-1}b+\frac{\eta}{h})'Z^i+\frac{1}{2\beta h^2}|h\sigma^{-1}b+\eta|^2\\
\qquad\qquad\quad+\frac{\rho^i}{\beta}-\frac{h}{\beta}(1-\ln h)-\frac{1}{\beta}\sum_{j=1}^\ell q^{ij}e^{-\beta(Y^j-Y^i)}\Big]dt+(Z^i)'dW, \\
Y^i_T=0, \quad i\in\cM.
\end{cases}
\end{align}
\begin{remark}
Because $\eta\in\BMO$ is not bounded, we cannot apply \cite[Theorem 2.4]{FHT} to solve \eqref{Yexp:ran}.
\end{remark}

Set $(P^i, \Lambda^i)=(\frac{Y^i}{h}, \frac{1}{h}Z^i-\frac{Y^i}{h^2}\eta)$ for all $i\in\cM$. Then $(P^i, \Lambda^i)_{i\in\cM}$ fulfills the following BSDE system:
\begin{align}\label{Pexp:ran}
\begin{cases}
dP^i=-\Big[-rP^i-(\sigma^{-1}b)'\Lambda^i+\frac{1}{2\beta h^3}|h\sigma^{-1}b+\eta|^2\\
 \qquad\qquad\quad+\frac{\rho^i}{\beta h}-\frac{1}{\beta}(1-\ln h)-\frac{1}{\beta h}\sum_{j=1}^\ell q^{ij}e^{-\beta h(P^j-P^i)}\Big]dt+(\Lambda^i)'dW, \\
P^i_T=0, \quad i\in\cM.
\end{cases}
\end{align}

\begin{theorem}\label{Th:exp:P}
Under Assumptions \ref{assu:r:ran} and \ref{assu:r:ran2}, \eqref{Pexp:ran} admits a unique solution $(P^i, \Lambda^i)_{i\in\cM}$ such that $(P^i, \Lambda^i)\in L^\infty_{\mathcal{F}^W}(0, T; \mathbb{R})\times \BMO$ for all $i\in\cM$.
\end{theorem}
\pf
We take the notations in the proof of Theorem \ref{Th:exp:h}.
Under the probability measure $\widetilde {\mathbb{P}}$,
\[
\widetilde W_t:=W_t+\int_0^t \sigma^{-1}b ds,
\]
 is an $n$-dimensional Brownian motion.
Notice that $\eta\in\BMO$, hence $\int_0^t \eta'd \widetilde W_s$ is a BMO martingale under $\widetilde {\mathbb{P}}$.

Let $k_1>0$ be a large constant such that, for all $i\in\cM$ and $t\in[0, T]$,
\begin{align}
 -r\leq k_1, \quad \frac{\rho^i}{\beta h}-\frac{1}{\beta}(1-\ln h)\geq -k_1,
\end{align}
and
\begin{align*}
&\widetilde\E\bigg[\int_t^T e^{-\int_t^s r_udu}\Big(\frac{1}{2\beta h^3}|h\sigma^{-1}b+\eta|^2+\frac{\rho^i}{\beta h} -\frac{1}{\beta}(1-\ln h)-\frac{1}{\beta h}q^{ii}\Big)ds\;\bigg|\;\mathcal{F}^W_{t}\bigg]\leq k_1.
\end{align*}
Denote $k:=\max\{e^{k_1 T}, k_1\}$.

Define a truncating function $\psi: \R\rightarrow\R$ as
\[
\psi(y):=\max\{-k, \min\{y, k\}\}.
\]
Consider the following truncated system of \eqref{Pexp:ran},
\begin{align}
\label{Ptrun}
\begin{cases}
dP^i=-\Big[-rP^i-(\sigma^{-1}b)'\Lambda^i+\frac{1}{2\beta h^3}|h\sigma^{-1}b+\eta|^2\\
\qquad\qquad\quad+\frac{\rho^i}{\beta h}-\frac{1}{\beta}(1-\ln h)-\frac{1}{\beta h}\sum_{j=1}^\ell q^{ij}e^{-\beta h(\psi(P^j)-\psi(P^i))}\Big]dt+(\Lambda^i)'dW, \\
P^i_T=0, \quad i\in\cM.
\end{cases}
\end{align}
The driver is Lipschitz continuous except for $\eta\in\BMO$ is unbounded. By \cite[Theorem 9.3.5]{CZ}, there is a unique solution $(P^{k, i}, \Lambda^{k, i})_{i\in\cM}$, such that $(P^{k, i}, \Lambda^{k, i})\in L^2_{\mathcal{F}^W}(0, T; \mathbb{R})\times L^2_{\mathcal{F}^W}(0, T; \mathbb{R}^n)$ for all $i\in\cM$.
If we can show $|P^{k, i}|\leq k$ and $\Lambda^{k, i}\in \BMO$, then $\psi(P^{k, i})=P^{k, i}$ and $(P^{k, i}, \Lambda^{k, i})_{i\in\cM}$ is a solution to \eqref{Pexp:ran}.

 We first show that $P^{k, i}\geq -k$.
It is direct to verify that
\begin{align*}
(\underline P, \underline \Lambda)=(1-e^{k_1(T-t)}, 0),
\end{align*}
is the unique solution to the following BSDE:
\begin{align}
\begin{cases}
d\underline P=-(k_1 \underline P-k_1)dt+\underline \Lambda'dW \\
\underline P_T=0.
\end{cases}
\end{align}
Notice that $\underline P\leq 0$ and
\begin{align*}
k_1\underline P-k_1&\leq -r\underline P-(\sigma^{-1}b)'\underline \Lambda+\frac{1}{2\beta h^3}|h\sigma^{-1}b+\eta|^2\\
&\qquad+\frac{\rho^i}{\beta h}-\frac{1}{\beta}(1-\ln h)-\frac{1}{\beta h}\sum_{j=1}^\ell q^{ij}e^{-\beta h(\psi(\underline P)-\psi(\underline P))},
\end{align*}
so by Lemma \ref{comparison} $$-k\leq -e^{k_1T}\leq\underline P_t\leq P^{k, i}_t.$$

We next show that $P^{k, i}\leq k$. From \eqref{Ptrun}, we obtain
\begin{align*}
e^{-\int_0^t r_sds}P^{k, i}_t&=\int_t^T e^{-\int_0^s r_udu}\Big[-(\sigma^{-1}b)'\Lambda^{k, i}+\frac{1}{2\beta h^3}|h\sigma^{-1}b+\eta|^2\\
&\qquad+\frac{\rho^i}{\beta h}-\frac{1}{\beta}(1-\ln h)-\frac{1}{\beta h}\sum_{j=1}^\ell q^{ij}e^{-\beta h(\psi(P^{k, j})-\psi(P^{k, i}))}\Big]ds-\int_t^T (\Lambda^{k, i})'dW\\
&=\int_t^T e^{-\int_0^s r_udu}\Big[\frac{1}{2\beta h^3}|h\sigma^{-1}b+\eta|^2+\frac{\rho^i}{\beta h}-\frac{1}{\beta}(1-\ln h)\\
&\qquad-\frac{1}{\beta h}\sum_{j=1}^\ell q^{ij}e^{-\beta h(\psi(P^{k, j})-\psi(P^{k, i}))}\Big]ds-\int_t^T (\Lambda^{k, i})'d\widetilde W.
\end{align*}
Taking expectation $\widetilde\E$, and noting $q^{ij}\geq 0$ for $i\neq j$, we have
\begin{align*}
P^{k, i}_t&=e^{\int_0^t r_sds}\widetilde\E\bigg[\int_t^T e^{-\int_0^s r_udu}\Big(\frac{1}{2\beta h^3}|h\sigma^{-1}b+\eta|^2+\frac{\rho^i}{\beta h}-\frac{1}{\beta}(1-\ln h)\\
&\qquad\qquad\qquad-\frac{1}{\beta h}q^{ii}-\frac{1}{\beta h}\sum_{j\neq i}^\ell q^{ij}e^{-\beta h(\psi(P^{k, j})-\psi(P^{k, i}))}\Big)ds\;\bigg|\;\mathcal{F}^W_{t}\bigg]\\
&\leq \widetilde\E^{i}\bigg[\int_t^T e^{-\int_t^s r_udu}\Big(\frac{1}{2\beta h^3}|h\sigma^{-1}b+\eta|^2+\frac{\rho^i}{\beta h}-\frac{1}{\beta}(1-\ln h)-\frac{1}{\beta h}q^{ii}\Big)ds\;\bigg|\;\mathcal{F}^W_{t}\bigg]\\
&\leq k_1\leq k.
\end{align*}
Therefore $|P^{k, i}|\leq k$. Whence we can drop the superscript $k$ in $(P^{k, i}, \Lambda^{k, i})_{i\in\cM}$.

Applying It\^{o}'s formula to $(P^i)^2$, we have for any stopping time $\tau\leq T$,
\begin{align*}
&\quad( P^i_{\tau})^2+\E\bigg[\int_{\tau}^T | \Lambda^i|^2ds\;\bigg|\;\mathcal{F}^W_{\tau}\bigg]\\
&=\E\bigg[\int_{\tau}^T 2P^i\Big(-rP^{i}-(\sigma^{-1}b)'\Lambda^{i}+\frac{1}{2\beta h^3}|h\sigma^{-1}b+\eta|^2\\
&\qquad\qquad\quad+\frac{\rho^i}{\beta h}-\frac{1}{\beta}(1-\ln h)-\frac{1}{\beta h}\sum_{j=1}^\ell q^{ij}e^{-\beta h(P^{j}-P^{i})}\Big)ds\;\bigg|\;\mathcal{F}^W_{\tau}\bigg]\\
&\leq \E\bigg[\int_{\tau}^T \Big(K+\frac{1}{2}|\Lambda^i|^2+K|\eta|^2\Big)ds\;\bigg|\;\mathcal{F}^W_{\tau}\bigg].
\end{align*}
From which, we obtain
\begin{align*}
\E\bigg[\int_{\tau}^T | \Lambda^i|^2ds\;\bigg|\;\mathcal{F}^W_{\tau}\bigg]\leq K\E\bigg[\int_{\tau}^T \Big(1+|\eta|^2\Big)ds\;\bigg|\;\mathcal{F}^W_{\tau}\bigg]\leq K,
\end{align*}
therefore $\Lambda^i\in\BMO$. The uniqueness comes from \cite[Lemma 3.4]{HSX}.
\eof

From Theorem \ref{Th:exp:h} and Theorem \ref{Th:exp:P}, we immediately have
\begin{theorem}
Under Assumptions \ref{assu:r:ran} and \ref{assu:r:ran2}, \eqref{Yexp:ran} admits a unique solution $(Y^i, Z^i)_{i\in\cM}$ such that $(Y^i, Z^i)\in L^\infty_{\mathcal{F}^W}(0, T; \mathbb{R})\times \BMO$ for all $i\in\cM$.
\end{theorem}

\begin{theorem}\label{Th:opt}
Suppose Assumptions \ref{assu:r:ran} and \ref{assu:r:ran2} hold. Let $(h, \eta)$ and
$(Y^i, Z^i)_{i\in\cM}$ be the unique solutions to \eqref{hexp:ran} and \eqref{Yexp:ran} respectively. Then the value function of the optimization problem \eqref{optmexp} is given by
$$V(x, i_0)=-e^{-\beta(h_0 x+Y_0^{i_0})},$$ and the optimal portfolio is
\begin{align}\label{optmpor:exp:r}
\hat\pi_t=-\frac{1}{\beta h^2_t} (\sigma'_t)^{-1}(\beta h_t\eta_t X_t+\beta h_tZ_t^{\alpha_t}-h_t\sigma_t^{-1}b_t-\eta_t),
\end{align}
the optimal consumption is
\begin{align}\label{optmcon:exp:r}
\hat c_t=h_t X_t+Y_t^{\alpha_t}-\frac{1}{\beta}\ln h_t.
\end{align}
\end{theorem}
\pf
For any $(\pi, c)\in\mathcal{U}$, applying It\^{o}'s formula to $-e^{-\int_0^t\rho_s^{\alpha_s}ds}e^{-\beta(h_tX_t+Y_t^{\alpha_t})}$,
we obtain 
\begin{align*} \E\Big[-e^{-\int_0^{\tau_n}\rho_s^{\alpha_s}ds}e^{-\beta(h_{\tau_n}X_{\tau_n}+Y_{\tau_n}^{\alpha_{\tau_n}})}\Big]
-\E\Big[\int_0^{\tau_n} e^{-\int_0^s\rho_s^{\alpha_u}du}e^{-\beta c_s}ds\Big]\leq -e^{-\beta(h_0x+Y_0^{i_0})},
\end{align*}
for a sequence of increasing stopping times $(\tau_n)_{n\in\mathbb{N}}$ satisfying $\lim_{n\rightarrow\infty}\tau_n=T$. 
Same as the deterministic case, the two expectations are both convergence as $n\to\infty$ so that
\begin{align} \label{superma}
\E\Big[-e^{-\int_0^{T}\rho_s^{\alpha_s}ds}e^{-\beta X_{T}}\Big]
-\E\Big[\int_0^{T} e^{-\int_0^s\rho_s^{\alpha_u}du}e^{-\beta c_s}ds\Big]\leq -e^{-\beta(h_0x+Y_0^{i_0})},
\end{align}
thanks to the terminal conditions $h_T=1$ and $Y_T=0$. Moreover, the inequality \eqref{superma} becomes an equality for some $(\hat\pi, \hat c)$ if and only if \eqref{optmpor:exp:r} and \eqref{optmcon:exp:r} hold.

Denote $\hat X$ by the corresponding wealth process under the strategy $(\hat\pi, \hat c)$. Then we have
\begin{align*}
\E\Big[-e^{-\int_0^{\tau_n}\rho_s^{\alpha_s}ds}e^{-\beta(h_{\tau_n}^{\alpha_{\tau_n}}\hat X_{\tau_n}+Y_{\tau_n}^{\alpha_{\tau_n}})}
-\int_0^{\tau_n} e^{-\int_0^s\rho_s^{\alpha_u}du}e^{-\beta \hat c_s}ds\Big]=-e^{-\beta(h_0x+Y_0^{i_0})}.
\end{align*}
Note that $\rho$, $Y^{i}$, $i\in\cM$, are bounded processes,
\begin{align*}
 \E\Big[e^{-\beta h_{\tau_n}^{\alpha_{\tau_n}}\hat X_{\tau_n}}
+\int_0^{\tau_n} e^{-\beta \hat c_s}ds\Big]
&\leq K \E\Big[e^{-\int_0^{\tau_n}\rho_s^{\alpha_s}ds}e^{-\beta(h_{\tau_n}^{\alpha_{\tau_n}}\hat X_{\tau_n}+Y_{\tau_n}^{\alpha_{\tau_n}})}
+\int_0^{\tau_n} e^{-\int_0^s\rho_s^{\alpha_u}du}e^{-\beta \hat c_s}ds\Big]\\
&=K e^{-\beta(h_0x+Y_0^{i_0})}.
\end{align*}
Applying Fatou's lemma to the above, we obtain
\begin{align}\label{expint}
 \E\Big[e^{-\beta \hat X_{T}}
+\int_0^{T} e^{-\beta \hat c_s}ds\Big]\leq K.
\end{align}

Set $\theta:=\frac{1}{h}\sigma^{-1}(hb+\sigma\eta)$, then $\theta\in \BMO$ and $\mathcal{E}(\int_0^t\theta_sdW)$ is a uniformly integrable martingale. By Girsanov's theorem,
\[
W^\mathbb{Q}_t:=W_t+\int_0^t\theta_sds
\]
is an $n$-dimensional Brownian motion under the measure $\mathbb{Q}$ defined by
\[
\frac{d\mathbb{Q}}{d\mathbb{P}}=\mathcal{E}\Big(-\int_0^T\theta_s dW_s\Big).
\]

We need to show that $(e^{-\beta h_t\hat X_t})_{0\leq t\leq T}$ belongs to class (D). Compared with Section 4 or \cite{CH}, a new term $\eta\in\BMO$ emerges (In Section 4 or \cite{CH}, $r$ is deterministic, hence $\eta=0$).
Applying It\^{o}'s formula to $h_t \hat X_t$,
\begin{align*}
h_t \hat X_t&=h_0x+\int_0^t\Big[-\frac{1}{\beta h^2}(\beta hZ-h\sigma^{-1}b-\eta)'\sigma^{-1}(hb+\sigma\eta)-h(Y-\frac{1}{\beta}\ln h)\Big]ds\\
&\qquad\quad-\int_0^t\frac{1}{\beta h}(\beta hZ-h\sigma^{-1}b-\eta)'dW\\
&=h_0x+\int_0^t\Big[-h(Y-\frac{1}{\beta}\ln h)\Big]ds-\int_0^t\frac{1}{\beta h}(\beta hZ-h\sigma^{-1}b-\eta)'dW^{\mathbb{Q}}
\end{align*}
Since $-\frac{1}{\beta h}(\beta hZ-h\sigma^{-1}b-\eta)\in L^{2, \;\mathrm{BMO}}_{\mathcal F}(0, T;\mathbb{R}^{n})$, the process $$L_t:=-\int_0^t\frac{1}{\beta h}(\beta hZ-h\sigma^{-1}b-\eta)'dW^{\mathbb{Q}}$$ is a BMO-martingale under $\mathbb{Q}$. Because $h$, $Y^{i}$, $i\in\cM$, are bounded, there exists a positive constant $K$ such that
\begin{align*}
\frac{1}{K}e^{-\beta L_t}\leq e^{-\beta h_t \hat X_t}\leq K e^{-\beta L_t}
\ \mbox{for all} \ t\in[0, T].
\end{align*}
Hence our problem reduces to show $(e^{-\beta L_t })_{0\leq t\leq T}$ belongs to class (D).
By the reverse H\"older inequality, there exists $p>1$ such that for all stopping times $\tau\leq T$,
\begin{align}\label{reverseholder}
\Big(\E\Big[\mathcal{E}\Big(-\int_{0}^T\theta dW\Big)^p\;\Big|\;\mathcal{F}_{\tau}\Big]\Big)^{\frac{1}{p}}\leq K \mathcal{E}\Big(-\int_{0}^{\tau}\theta dW\Big).
\end{align}
Let $q$ be the conjugate exponent of $p$ i.e. $\frac{1}{p}+\frac{1}{q}=1$. Since the function $x\mapsto e^{-\frac{\beta}{q}x}$ is convex, by Jensen's inequality, we have for any stopping time $\tau\leq T$,
\begin{align*}
e^{-\frac{\beta}{q}L_{\tau}}\leq \E^{\mathbb{Q}}\Big[e^{-\frac{\beta}{q}L_T}\;\Big|\;\mathcal{F}_{\tau}\Big].
\end{align*}
It follows
\begin{align*}
 e^{-\beta L_{\tau}}
&\leq \Big(\E^{\mathbb{Q}}\Big[e^{-\frac{\beta}{q}L_T}\;\Big|\;\mathcal{F}_{\tau}\Big]\Big)^q\\
&= \Big(\E\Big[e^{-\frac{\beta}{q}L_T}\mathcal{E}\Big(-\int_0^T \theta dW \Big)\;\Big|\;\mathcal{F}_{\tau}\Big] \Big)^q\mathcal{E}\Big(-\int_0^{\tau}\theta dW \Big)^{-q}\\
&\leq \E\Big[e^{- \beta L_T}\;\Big|\;\mathcal{F}_{\tau}\Big] \E\Big[\mathcal{E}\Big(-\int_0^T \theta dW \Big)^p\;\Big|\;\mathcal{F}_{\tau}\Big]^{\frac{q}{p}}\mathcal{E}\Big(-\int_0^{\tau} \theta dW \Big)^{-q}\\
&\leq K \E[e^{-\beta L_T}\;\Big|\;\mathcal{F}_{\tau}],
\end{align*}
where we use H\"older's inequality in the second inequality, \eqref{reverseholder} in the last inequality. This shows that $(e^{-\beta h_t\hat X_t})_{0\leq t\leq T}$ belongs to class (D) noting \eqref{expint}, hence the admissibility of $(\hat\pi, \hat c)$.
\eof

\section{Concluding remarks}\label{conclude}
In this paper, we studied optimal consumption-investment problems in a market with regime switching and random coefficients with coupled constraints. Explicit solutions in term of the solutions of BSDE systems are given for the cases of power, logarithmic and exponential utilities. Some of these BSDE systems are new in the literature. Their solvability consists the main mathematical contributions and is interesting in its own from the point of view of BSDE theory. Extensions in other directions can be interesting as well; for instance: (1) The problem in infinite time horizon (lifetime). Actually the discount factor $\rho$ was introduced in problems \eqref{optm} and \eqref{optmexp} for this possible extension. (2) The problem in an incomplete market, i.e. $m<n$ in subsection 5.2 or when $r$, $\mu$, or $\sigma$ depends on the regime $i$.

\end{document}